\documentclass[12pt,a4paper]{article}
\usepackage{latexsym,epsfig,fullpage,amsfonts,amssymb,amsmath,amscd,epic,graphics,rotating,theorem}
\theoremstyle{plain}
\newtheorem{lemma}{Lemma}[section]
\newtheorem{proposition}[lemma]{Proposition}
\newtheorem{remark}[lemma]{Remark}

\newtheorem{theorem}[lemma]{Theorem}
\newtheorem{definition}[lemma]{Definition}
\newtheorem{corollary}[lemma]{Corollary}
{\theorembodyfont{\rmfamily}
\begin{document}
\newcommand{\pperp}{\hbox{$\perp\hskip-6pt\perp$}}
\newcommand{\ssim}{\hbox{$\hskip-2pt\sim$}}
\newcommand{\N}{{\mathbb N}}
\newcommand{\A}{{\mathbb A}}
\newcommand{\Z}{{\mathbb Z}}
\newcommand{\R}{{\mathbb R}}
\newcommand{\C}{{\mathbb C}}
\newcommand{\Q}{{\mathbb Q}}
\newcommand{\Id}{{\operatorname{Id}}}
\newcommand{\real}{{\operatorname{Re}}}
\newcommand{\diag}{{\operatorname{diag}}}
\newcommand{\Span}{{\operatorname{Span}}}
\newcommand{\Fix}{{\operatorname{Fix}}}
\newcommand{\sign}{{\operatorname{sign}}}
\newcommand{\Tors}{{\operatorname{Tors}}}
\newcommand{\oi}{{\overline i}}
\newcommand{\oj}{{\overline j}}
\newcommand{\ob}{{\overline b}}
\newcommand{\os}{{\overline s}}
\newcommand{\oa}{{\overline a}}
\newcommand{\oy}{{\overline y}}
\newcommand{\ow}{{\overline w}}
\newcommand{\ot}{{\overline t}}
\newcommand{\oz}{{\overline z}}
\newcommand{\eps}{{\varepsilon}}
\newcommand{\proofend}{$\Box$\bigskip}
\newcommand{\Int}{{\operatorname{Int}}}
\newcommand{\pr}{{\operatorname{pr}}}
\newcommand{\grad}{{\operatorname{grad}}}
\newcommand{\rk}{{\operatorname{rk}}}
\newcommand{\im}{{\operatorname{Im}}}
\newcommand{\sk}{{\operatorname{sk}}}
\newcommand{\const}{{\operatorname{const}}}
\newcommand{\Sing}{{\operatorname{Sing}}}
\newcommand{\conj}{{\operatorname{Conj}}}
\newcommand{\Cl}{{\operatorname{Cl}}}
\newcommand{\discr}{{\operatorname{discr}}}
\newcommand{\Tor}{{\operatorname{Tor}}}
\newcommand{\defect}{{\operatorname{def}}}
\newcommand{\tmu}{{\C\mu}}
\newcommand{\ov}{{\overline v}}
\newcommand{\ox}{{\overline{x}}}
\newcommand{\tet}{{\theta}}
\newcommand{\Del}{{\Delta}}
\newcommand{\bet}{{\beta}}
\newcommand{\kap}{{\kappa}}
\newcommand{\del}{{\delta}}
\newcommand{\sig}{{\sigma}}
\newcommand{\alp}{{\alpha}}
\newcommand{\Sig}{{\Sigma}}
\newcommand{\Gam}{{\Gamma}}
\newcommand{\gam}{{\gamma}}
\newcommand{\Lam}{{\Lambda}}
\newcommand{\lam}{{\lambda}}
\newcommand{\nek}{{,...,}}
\newcommand{\inote}{\newline \centerline {It: \it{CHANGED}} \newline}
\title{VIRO THEOREM AND
TOPOLOGY OF REAL AND COMPLEX COMBINATORIAL
HYPERSURFACES\footnote{2000 Mathematics Subject Classification:
Primary 14M25, 14P25. Secondary 52B20, 52B70}}
\author{Ilia Itenberg\\
CNRS,
Institut de Recherche Math\'ematique de Rennes \\
Campus de Beaulieu \\
35042 Rennes Cedex, France\\
E-mail: itenberg@maths.univ-rennes1.fr
\and Eugenii Shustin
\thanks{A part of the present
work was done during the stay of the second author at the Fields
Institute, Toronto, and at the NSF Science and Technology Research
Center for the Computation and Visualization of Geometric Structures,
funded by NSF/DMS89-20161.
The work was completed  during the stay of both authors at
Max-Planck-Institut f\"ur Mathematik. The authors thank these
funds and institutions for hospitality and financial support.}
\\ Tel Aviv University\\
School of Mathematical Sciences\\
Ramat Aviv, 69978 Tel Aviv, Israel\\
E-mail: shustin@post.tau.ac.il}
\date{}
\maketitle
\begin{abstract}
We introduce a class of combinatorial hypersurfaces in the
complex projective space. They are submanifolds of codimension~$2$
in $\C P^n$
and are topologically "glued" out of algebraic hypersurfaces
in $(\C^*)^n$. Our construction can be viewed as a version
of the Viro
gluing theorem, relating topology of algebraic hypersurfaces to the
combinatorics of subdivisions of convex lattice polytopes.
If a subdivision is convex, then according
to the Viro theorem a combinatorial hypersurface is
isotopic to an algebraic one. We study combinatorial
hypersurfaces resulting from non-convex subdivisions of convex polytopes,
show that they are almost complex varieties, and in the real case,
they satisfy the same topological restrictions
(congruences, inequalities etc.)
as real algebraic hypersurfaces.
\end{abstract}

\section*{Introduction}

Topology of real algebraic varieties, brought to
the wide mathematical audience by D.~Hilbert in his
16th problem, has been studied from two sides,
looking for restrictions to geometric and topological
properties of real algebraic varieties, and constructing
varieties with prescribed properties (a lot of material, but,
certainly not all, can be found in the surveys \cite{De-Kh,Gu,Vi4,Wi}).
Analyzing a gap between restrictions and constructions, O.~Viro
\cite{Vi4} suggested a concept of ``flexible" curves, smooth surfaces
in $\C P^2$
having certain topological properties of
real algebraic curves, and asked for
their classification in comparison with the classification of
real algebraic curves. We mention a recent active study
of real pseudo-holomorphic
curves \cite{FO,O,OS,W} as a further development of this program.

In the present paper we suggest another source for flexible
curves, or, more generally, flexible varieties which we call {\it
combinatorial hypersurfaces}. These varieties are codimension $2$
submanifolds in the complex projective space, they are
conjugation invariant and their real parts are codimension $1$
submanifolds in the real projective space. Our approach is based
on the Viro construction of real algebraic varieties with
prescribed topology \cite{Vi1, Vi2, Vi3, Vi7} (see also \cite{GKZ},
11.5, \cite{IV, R}), which relates topology of real algebraic
varieties to the combinatorics of Newton polytopes and their {\it
convex} lattice subdivisions.
The Viro gluing theorem states that the result
of the gluing procedure is isotopic to a real algebraic
variety under the condition of convexity of the subdivision
used in the construction.
We remove the convexity condition and show that the Viro
construction modified in this way still can be performed and
produces combinatorial hypersurfaces. Then we show that
combinatorial hypersurfaces obey almost all known topological
restrictions to real algebraic varieties. These results can be
viewed as the first steps in the study of the
following questions: how far are combinatorial
hypersurfaces from the algebraic ones, and how crucial is the
convexity condition in the Viro theorem~?

We preface the main material with illustrating examples
and more detailed statement of the problem and our results.

Consider an example of the Viro construction.
Let $T_d\subset\R^2$, $d\in\N$, be the triangle with
vertices $(0,0)$, $(0,d)$, $(d,0)$, $$\tau:\quad
T_d=\Del_1\cup...\cup\Del_N$$ a triangulation with the set of
vertices $V\subset\Z^2$, and $\sig:V\to\{\pm 1\}$ any function.
Out of this combinatorial data we construct a piecewise-linear plane
curve. Denote by $T_d^{(1)}$, $T_d^{(2)}$ and $T_d^{(3)}$ the
copies of $T_d$ under reflections with respect to the coordinate
axes and the origin. Take in $T_d^{(1)}$, $T_d^{(2)}$ and
$T_d^{(3)}$ triangulations symmetric to~$\tau$, and define $\sig$
at the vertices of new triangulations by
$$\sig(\eps_1i,\eps_2j)=\eps_1^i\eps_2^j\sig(i,j),\quad (i,j)\in
V,\quad \eps_1,\eps_2=\pm 1.$$ Now in any triangle of the
triangulation of $T_d\cup T_d^{(1)}\cup T_d^{(2)}\cup T_d^{(3)}$,
having vertices with different values of $\sig$, we draw the
midline separating the vertices with different signs. The union
$C(\tau,\sig)$ of all these midlines is a broken line homeomorphic
to a disjoint union of circles and segments. Introduce the
following maps: $$\Phi:T_d\cup T_d^{(1)}\cup T_d^{(2)}\cup
T_d^{(3)}\to\R P^2,\quad \Psi:\Int(T_d\cup T_d^{(1)}\cup
T_d^{(2)}\cup T_d^{(3)})\to\R^2,$$ where $\Phi$ is continuous
onto, identifying antipodal points on $\partial(T_d\cup
T_d^{(1)}\cup T_d^{(2)}\cup T_d^{(3)})$, and $\Psi$ is a
homeomorphism. We call the curves
$\Phi(C(\tau,\sig))\subset\R P^2$
and $\Psi(C(\tau,\sig) \cap \Int(T_d\cup T_d^{(1)}\cup
T_d^{(2)}\cup T_d^{(3)})) \subset\R^2$
{\it projective} and {\it affine T-curves of degree}~$d$,
respectively.

The Viro theorem states
that a projective (resp., affine) T-curve
of degree~$d$ is isotopic in~$\R P^2$ (resp., in~$\R^2$) to a
nonsingular algebraic projective (resp., affine) curve of
degree~$d$, providing that the triangulation~$\tau$ is {\it
convex}, {\it i.e.}, there exists a convex piecewise-linear
function $\nu: T_d\to\R$, whose linearity domains are $\Del_1,
\ldots , \Del_N$ (sometimes such triangulations are called {\it
regular} or {\it coherent}; see~\cite{Z} and~\cite{GKZ}). The Viro
theorem, in fact, endows the combinatorial broken line
$C(\tau,\sig)$ with a rich structure, which yields a number of
restrictions on the topology of $C(\tau,\sig)$ (see an account of
known results in \cite{Ro, Vi5, Vi4, Wi}).

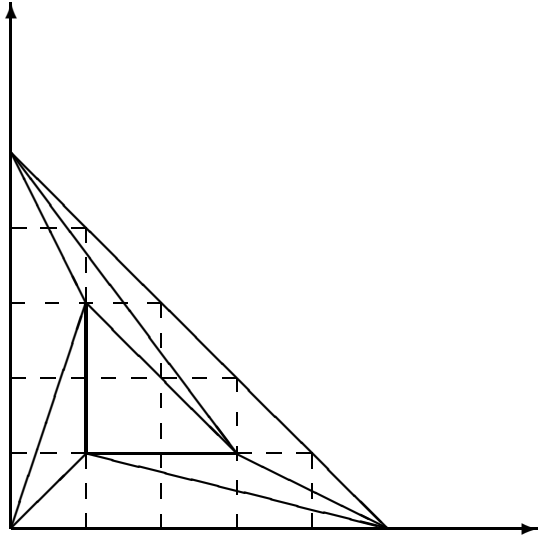
\begin{figure}
\setlength{\unitlength}{1cm}
\begin{picture}(8,8)(0,0)
\thicklines
\put(3,1){\vector(1,0){7}}
\put(3,1){\vector(0,1){7}}
\put(3,6){\line(1,-1){5}}
\put(3,6){\line(1,-2){1}}
\put(3,1){\line(1,3){1}}
\put(3,1){\line(1,1){1}}
\put(4,2){\line(4,-1){4}}
\put(4,2){\line(1,0){2}}
\put(4,2){\line(0,1){2}}
\put(4,4){\line(1,-1){2}}
\put(3,6){\line(3,-4){3}}
\put(6,2){\line(2,-1){2}}
\thinlines
\dashline{0.2}(3,2)(4,2)
\dashline{0.2}(6,2)(7,2)
\dashline{0.2}(3,3)(6,3)
\dashline{0.2}(3,4)(5,4)
\dashline{0.2}(3,5)(4,5)
\dashline{0.2}(4,1)(4,2)
\dashline{0.2}(4,4)(4,5)
\dashline{0.2}(5,1)(5,4)
\dashline{0.2}(6,1)(6,3)
\dashline{0.2}(7,1)(7,2)
\end{picture}
\caption{A non-convex triangulation}\label{f2}
\end{figure}

On the other hand, there exist {\it non-convex} triangulations.
The simplest example is shown in Figure~\ref{f2} (see, for
instance, \cite{CH}). Moreover, no efficient criterion for the
convexity of a triangulation is known.
There are examples of T-curves beyond the range of known
algebraic curves~\cite{Sa}, and there is some similarity between
T-curves and algebraic curves: up to degree~$6$ any projective
T-curve is isotopic to an algebraic one of the same
degree~\cite{dL1}, and {\it vice versa}, any nonsingular algebraic
curve of degree at most~$6$ in $\R P^2$ is isotopic to a
projective T-curve of the same degree. T-curves satisfy some
consequences of the B\'ezout theorem~\cite{dL1}, the Harnack
inequality \cite{I1, Ha},
and the complex orientation formula \cite{Pa}.
Maximal T-curves satisfy the Ragsdale type inequality \cite{Ha},
which is not proved for maximal real algebraic curves.
It is natural to ask
whether any T-curve
is isotopic to an algebraic curve of the same degree, and if not,
how far T-curves may differ from algebraic ones.

In the general Viro construction,
which applies to any dimension, one starts with
a convex lattice subdivision
$P = \Delta_1 \cup \ldots \cup \Delta_N$
of an $n$-dimensional convex lattice polytope~$P$ in $(\R_+)^n$
(where $\R_+ = \{x \in \R, x \geq 0\}$)
and an appropriate collection of polynomials $F_1$, $\ldots$, $F_N$.
Then the Viro construction
produces an algebraic hypersurface~$A$
in the toric variety~$\Tor_\C(P)$
associated with~$P$. The Newton polytope of~$A$ is~$P$,
and the topology of the complex part of~$A$
(and of the real part of~$A$ if the polynomials $F_1$, $\ldots$, $F_N$
are real) is described in terms of topology of zero point sets
of $F_1$, $\ldots$, $F_N$. Namely, the complex (resp., real)
part of~$A$ is,
in a sense, ``glued'' out of the complex
(resp., real) zero sets of $F_1$, $\ldots$, $F_N$.

The zero sets of the initial polynomials in the Viro construction
can be ``glued'' even if the subdivision is not convex
(as it was done above in the case of T-curves).
However, in this situation the Viro theorem does not guarantee that
the result carries an algebraic structure.
In the present paper we study the following question:
what can be said about the result of the ``gluing''
in the case of
{\it non-convex} subdivisions of Newton polytopes~?

Let us restrict ourselves to the situation when $P$ is the simplex
$T_d^n$ in $\R^n$ with vertices $(0,0,...,0)$, $(d,0,...,0)$, ...,
$(0,...,0,d)$. Note that the toric variety $\Tor_\C(T_d^n)$ associated
with $T_d^n$ is $\C P^n$. We show that the result of ``gluing'' of
the complex zero sets of the polynomials $F_1$, $\ldots$, $F_N$ is
a codimension~$2$ piecewise-smooth submanifold in $\C P^n$. We
call this piecewise-smooth submanifold a {\it (complex)
combinatorial hypersurface} (briefly, {\it C-hypersurface}) {\it
of degree} $d$ {\it in} $\C P^n$. If the initial polynomials are
real, then the resulting C-hypersurface is invariant under the
complex conjugation in $\C P^n$. In this case, the C-hypersurface
is called {\it real}. The real point set of a real C-hypersurface
is a piecewise-analytic hypersurface in $\R P^n$ and is the result
of ``gluing'' of the real zero sets of the polynomials $F_1$,
$\ldots$, $F_N$.

It is natural to compare the class of C-hypersurfaces
of degree~$d$ in $\C P^n$ with the class of algebraic hypersurfaces
of the same degree. The following questions arise.
\begin{enumerate}
\item {\it Is any C-hypersurface of degree~$d$ in~$\C P^n$
homeomorphic (or isotopic) to an algebraic hypersurface of the
same degree~?}
\item {\it Is the real part of any real C-hypersurface of degree~$d$
in~$\C P^n$ homeomorphic (or isotopic in $\R P^n$) to the real
part of a real algebraic hypersurface of the same degree~?}
\item {\it Is any real C-hypersurface of degree~$d$ in~$\C P^n$
equivariantly isotopic to a real algebraic hypersurface of the
same degree~?}
\end{enumerate}

In the present work we try to do the first steps to answering
these questions. We show that any C-hypersurface of degree~$d$ in
$\C P^n$ can be smoothed and carries an almost complex structure.
We prove that C-hypersurfaces of degree~$d$ in $\C P^n$ share a
lot of topological properties with algebraic hypersurfaces of
degree~$d$ in $\C P^n$. In particular, the answer to the first
question is positive in the case $n = 2$. Moreover, any real
C-curve in $\C P^2$ can be smoothed out into a {\it flexible}
curve in the sense of~\cite{Vi4}. As a corollary, we obtain that
arbitrary T-curves satisfy all {\it topological} restrictions
known for real algebraic curves. We also prove that any C-surface
of degree~$d$ in $\C P^3$ is homeomorphic to an algebraic surface
of the same degree.

The material of section 1, where we describe the construction of
C-hypersurfaces, is basically contained in \cite{Vi1}.
We provide here all details, since our setting is a little bit different.
In section 2 we study the
topological properties of (complex) C-hypersurfaces of degree~$d$
in~$\C P^n$ and of double coverings of~$\C P^n$ ramified along
C-hypersurfaces. In particular, we reprove the projective complex
version of Viro's theorem (see~\cite{Vi1}) which claims that
a C-hypersurface (resp., real C-hypersurface) of degree~$d$ in $\C P^n$
is isotopic (resp., equivariantly isotopic)
to an algebraic
hypersurface of the same degree, provided the subdivision is convex.
In the case of
arbitrary (not necessarily convex) subdivisions, we prove the
following statements.
\begin{itemize}
\item A C-hypersurface of degree~$d$ in~$\C P^n$ is an orientable manifold,
homologous to an algebraic
hypersurface of degree $d$ in $\C P^n$.
\item Any C-hypersurface~$M$ in~$\C P^n$ is
isotopic to a close smooth
manifold~$M_{sm}$ of codimension 2 in $\C P^n$.
If~$M$ is real, the isotopy can be made equivariant.
\item Given a (real) C-hypersurface $M$,
the tangent bundle
to its smoothing $M_{sm}$
is (equivariantly) isotopic to a
(equivariant) bundle of complex hyperplanes.
In particular, $M_{sm}$ possesses an (equivariant) almost complex structure.
\item A C-hypersurface in~$\C P^n$ is simply connected if $n>2$.
\item For a C-hypersurface~$M$ of degree~$d$ in~$\C P^n$, one has
$\pi_1(\C P^n\backslash M)=\Z/d\Z$ if $n\ge 2$.
\item Let $M \subset \C P^n$ be a C-hypersurface and
$M'\subset \C P^n$
a nonsingular algebraic hypersurface, both of degree
$d$. Then $H_*(M)$ and $H_*(M')$ are isomorphic as graded
groups.
\item Let~$n$ be a positive odd number, $M$ a C-hypersurface
of degree~$d$
in~$\C P^n$ and
$M'$ a nonsingular algebraic hypersurface of degree~$d$ in~$\C P^n$.
Then the lattices $(H_{n - 1}(M), B_M)$
and $(H_{n - 1}(M'), B_{M'})$
(where
$B_M: H_{n - 1}(M) \times H_{n - 1}(M) \to \Z$ and
$B_{M'}: H_{n - 1}(M') \times H_{n - 1}(M') \to \Z$ are
the intersection forms of~$M$ and~$M'$, respectively)
are isomorphic.
\item A C-surface of degree~$d$ in $\C P^3$ is homeomorphic to
a nonsingular algebraic surface of degree~$d$ in $\C P^3$.
\item Any C-curve in~$\C P^2$ is isotopic to an algebraic curve
of the same degree.
\end{itemize}
Section 3 is devoted to
topology of real C-hypersurfaces in~$\C P^n$.
We prove for them the generalized Harnack inequality and
the Gudkov-Rokhlin and Gudkov-Krahnov-Kharlamov congruences.
For real C-surfaces in~$\C P^3$, we also prove inequalities
similar to the Comessatti ones.

\begin{remark}
An interesting question concerns
the existence of a Hodge-like decomposition
in the cohomology of C-hypersurfaces.
A closely related question is
whether the analog of Petrovsky-Oleinik
inequalities is true
for real C-hypersurfaces in $\C P^n$~?
\end{remark}

\section{Construction of C-hypersurfaces}

\subsection{Notations and definitions}
Further on the term {\it polytope} ({\it polygon})
means a convex lattice polytope (polygon)
in the nonnegative orthant $\R^n_+$ of $\R^n$, $n\ge 2$.

Given a polynomial
$$F=\sum_{i_1,...,i_n}A_{i_1...i_n}z_1^{i_1}\cdot... \cdot
z_n^{i_n}\ ,$$ by $\Del(F)$ we denote its Newton polytope, {\it
i.e.}, the convex hall of the set $$\{(i_1,...,i_n)\in\R^n\ :\
A_{i_1...i_n}\ne 0\}\ .$$ The truncation of $F$ on a face $\del$
of $\Del(F)$ is the polynomial
$$F^{\del}=\sum_{(i_1,...,i_n)\in\del}A_{i_1...i_n}z_1^{i_1}\cdot...
\cdot z_n^{i_n}\ .$$ A polynomial $F\in\C[z_1,...,z_n]$ is called
{\it non-degenerate}, if $F$ and any truncation $F^{\del}$ on a
proper face $\del$ of $\Del(F)$ has a nonsingular zero set in
$(\C^*)^n$ (cf. \cite{Vi1}).

\subsection{Extension of the moment map}\label{s3}
Let $\Del\subset\R^n$ be a polytope, and let
$\mu_{\Del}:(\R^*_+)^n\to I(\Del)$
be the moment map (see \cite{A1, A2, Fu}, \cite{GKZ}, 6.1),
where $\R^*_+ = \{x \in \R, x > 0\}$ and $I(\Del)$ is the complement
in $\Del$ of the union of all its proper faces,
\begin{equation}
\mu_{\Del}(x_1,...,x_n)=
\frac{\sum_{(i_1,...,i_n)\in\Del}x_1^{i_1}...x_n^{i_n}
\cdot(i_1,...,i_n)}{\sum_{(i_1,...,i_n)\in\Del}x_1^{i_1}...
x_n^{i_n}}\ .\label{e1}
\end{equation}
Split the complex torus in the product $(\C^*)^n=(\R^*_+)^n
\times(S^1)^n$: $$(z_1,...,z_n)\in(\C^*)^n\ \mapsto \
(|z_1|,...,|z_n|)\in(\R^*_+)^n,\ \left(\frac{z_1}{|z_1|},...,
\frac{z_n}{|z_n|}\right)\in(S^1)^n\ .$$
 Note that the inverse map
$(\R^*_+)^n \times(S^1)^n\to(\C^*)^n$ naturally extends to a
surjection $\theta:\R^n_+\times(S^1)^n\to\C^n$. Put
 $$\C I(\Del) =\theta(I(\Del)\times(S^1)^n)\subset\C^n,
\quad \C\Del=\theta(\Del\times(S^1)^n)\subset\C^n\ .$$

\begin{proposition}\label{p3}
The complexification~$\C\Del$ of~$\Del$
is a singular PL-manifold with boundary.
The singular set of~$\C\Del$ is
the union of $\C\del$ over all faces $\del\subset\Del$,
which are intersections of $\Del$ with coordinate planes of
dimension $>n-\dim\Del+\dim\del$.
The real part $\R\Del$ of $\C\Del$
is the union of $\Del$ with all its symmetric copies
with respect to the coordinate hyperplanes.
\end{proposition}

{\bf Proof.} Straightforward. \proofend

Define {\it the extended moment map} $\tmu_{\Del}:(\C^*)^n\to
\C I(\Del)$ by
\begin{eqnarray}
&\tmu_{\Del}(x_1v_1,...,x_nv_n)=\theta(\mu_{\Del}(x_1,...,x_n)\ ,
(v_1,...,v_n)), \label{e2}\\
&\quad
(x_1,...,x_n)\in(\R^*_+)^n,\ (v_1,...,v_n)\in(S^1)^n\ ,
\nonumber \\
&\theta(\mu_{\Del}(x_1,...,x_n),
(v_1,...,v_n))\in\theta(I(\Del)\times(S^1)^n)=\C I(\Del)\ .
\nonumber
\end{eqnarray}
As an easy consequence of classical results we obtain
the following statement.

\begin{proposition}\label{p4}
The map $\tmu_{\Del}$ is surjective and
commutes with the
complex conjugation $\conj$. It is a diffeomorphism when
$\dim\Del=n$. The real part of $\C I(\Del)$ is the image of
$(\R^*)^n$.
\end{proposition}

Below we use
the following interaction of the extended moment map with
the action of $\Z^n$ in $\R^n$
$$\oa=(a_1,...,a_n)\in\Z^n,\ \ot=(t_1,...,t_n)\in\R^n\quad
\mapsto\quad\oa+\ot=(a_1+t_1,...,a_n+t_n)\in\R^n\ ,$$
and the actions of $SL(n,\Z)$ in $\R^n$ and $(\R^*)^n$
$$A\in SL(n,\Z),\ \ox\in\R^n\quad
\mapsto\quad A\ox\in\R^n\ ,$$
$$A\in SL(n,\Z),\ \ox=(x_1,...,x_n)\in(\R^*)^n
\ \mapsto\ \ox^A=\left(\prod_{i=1}^nx_i^{a_{i1}},...,
\prod_{i=1}^nx_i^{a_{in}}\right)\ .$$

\begin{proposition}\label{p7}
Let $\Del\subset\R^n$ be a polygon.
If $\oa=(a_1,...,a_n)
\in\Z^n$, $A\in SL(n,\Z)$, and
$\oa+A\Del$ lies in the nonnegative orthant $\R^n_+$, then
\begin{equation}
\tmu_{\oa+A\Del}(\oz)=(\oa+A\mu_\Del(\ox^A),\ov)\ ,
\label{e12}
\end{equation}
where $\oz=(\ox,\ov)\in(\R^*_+)^n\times(S^1)^n$.
\end{proposition}

{\bf Proof.} The result follows from the definition
$$\tmu_{\oa+A\Del}(\oz)=(\mu_{\oa+A\Del}(\ox),\ov)\ ,$$
and the classically known relation
$$\mu_{\oa+A\Del}(\ox)=\oa+A\mu_{\Del}(\ox^A)\quad
\mbox{\proofend}$$

\subsection{Real and complex chart of a polynomial}
Let $F\in\C[z_1,...,z_n]$ be a non-degenerate polynomial with
Newton polytope $\Del(F)=\Del$. The closure $\C
Ch(F)\subset\C\Del$ of the set $\tmu_{\Del}(\{F=0\}\cap(\C^*)^n)$
is called the {\it (complex) chart} of the polynomial~$F$. If~$F$
is real then $\R Ch(F)= \C Ch(F)\cap\R\Del$ is called the {\it
real chart} of~$F$.

This definition is a
key ingredient of the
Viro construction (see \cite{Vi1, R},
cf. \cite{Sh}).

\begin{proposition}\label{l1}
Suppose that $\Del\subset (\R^*_+)^n$.
Then the set $\C Ch(F)$ is a PL-submanifold
in $\C\Del$
of codimension 2 with boundary
$\partial\C Ch(F)=
\C Ch(F)\cap\partial\C\Del$.
It is smooth in $\C I(\Del)$, and,
for any proper face $\del$ of $\Del$ of positive dimension,
\begin{equation}
\C Ch(F)\cap\C\del=\C Ch(F^{\del})\ .\label{e13}
\end{equation}
If~$F$ is real then
$\C Ch(F)$ is invariant with respect to $\conj$, and $\R Ch(F)$ is a
PL-submanifold in~$\R\Del$ of codimension 1
with boundary
$\partial\R Ch(F)=
\R Ch(F)\cap\partial\R\Del$.
\end{proposition}

\begin{remark}
If $\Del$ intersects with coordinate hyperplanes, then the
statement of Proposition \ref{l1} holds true when substituting
$\C\Del\backslash\Sing(\C\Del)$ and $\R\Del\backslash\Sing(\C\Del)$
for $\C\Del$ and $\R\Del$, respectively, but we will not use this
below.
\end{remark}

{\bf Proof of Proposition \ref{l1}.}

(i) If $\dim\Del=n$ then by
Proposition~\ref{p4}, $\C Ch(F)\cap\C I(\Del)$ is a smooth manifold
of codimension 2.

So, we have to consider only the behavior
of $\C Ch(F)$ on $\partial\C\Del$.

(ii) Assume that $\dim\Del<n$ and $A\in SL(n,\Z)$ is
such that $A\Del\subset(\R^*_+)^n$ is parallel to the
coordinate hyperplane $i_n=0$. By (\ref{e12}) there is a
diffeomorphism which takes the pair $(\C\Del,\C Ch(F(\oz)))$ onto
the pair $(\C(A\Del),\C Ch(F(\oz^A)))$. Now, if $i_n\big|_\Del=a$,
then $\Del'=\Del-(0,...,0,a)\subset\{i_n=0\}$,
$F(\oz)=z_n^aF'(\oz)$, and by formula (\ref{e12})
the pair $(\C\Del,\C Ch(F))$ naturally splits into the product
$(C\Del',\C Ch(F'))\times S^1$, thereby one reduces the dimension of
the ambient space.

(iii) To prove (\ref{e13}) (for
$\dim\Del=n$), we note that
any point $\oz^0\in\C Ch(F)\cap
\partial\C\Del$ is a limit of
$\tmu_\Del(\gam(t))$ as $t>0$, $t\to 0$, where
$$\gam(t)=(\lam_1t^{k_1}+O(t^{k_1+1}),...,
\lam_nt^{k_n}+O(t^{k_n+1}))$$
is a curve lying in $\{F=0\}\cap(\C^*)^n$ with
$\lam_1...\lam_n\ne 0$, $(k_1,...,k_n)\in\Z^n\backslash\{0\}$.
Let $\del$ be the maximal (proper) face of $\Del$, where
the function $(i_1,...,i_n)\in\Del\mapsto i_1k_1 + ... + i_nk_n$
achieves its minimum (equal to $l$). Then
$$0=F(\gam(t))=F^\del(\lam_1,...,\lam_n)t^l+O(t^{l+1})\quad
\Longrightarrow\quad
F^\del(\lam_1,...,\lam_n)=0\ ,$$
$\gam(t)=(\ox,\ov)\in
(\R^*_+)^n\times(S^1)^n$, where
$$\ox=(|\lam_1|t^{k_1}+O(t^{k_1+1}),...,|\lam_n|t^{k_n}+O(t^{k_n+1})),
\ \ov=\left(\frac{\lam_1}{|\lam_1|}+O(t),...,
\frac{\lam_n}{|\lam_n|}+O(t)\right).$$
So, one derives
$$\tmu_\Del(\gam(t))=(\mu_\Del(\ox),\ov)=
\left(\frac{\sum_{(i_1,...,i_n)\in\Del}(i_1,...,i_n)x_1^{i_1}...
x_n^{i_n}}{\sum_{(i_1,...,i_n)\in\Del}x_1^{i_1}...
x_n^{i_n}},\ \ov\right)$$
$$=\left(\frac{t^l\cdot\sum_{(i_1,...,i_n)\in\del}(i_1,...,i_n)
|\lam_1|^{i_1}...|\lam_n|^{i_n}+O(t^{l+1})}{t^l\cdot
\sum_{(i_1,...,i_n)\in\del}|\lam_1|^{i_1}...
|\lam_n|^{i_n}+O(t^{l+1})},\
\ov\right)\quad\stackrel{t\to 0}{\longrightarrow}\quad
\tmu_\del(\lam_1,...,\lam_n).$$
Hence $\C Ch(F)\cap\C\del\subset\C Ch(F^{\del})$.

Now let $(\lam_1,...,\lam_n)\in(\C^*)^n$,
$F^\del(\lam_1,...,\lam_n)=0$,
where $\del$ is a proper face of $\Del$,
and let $(k_1,...,k_n)\in k(\del)$,
where $k(\del)\subset\Z^n\backslash\{0\}$ consists of vectors
$(k_1,...,k_n)$ such that $\del$ is the maximal face of
$\Del$ on which the function $(i_1,...,i_n)\in\Del
\mapsto i_1k_1+...+i_nk_n$ achieves its minimum.
Then $F^\del(\lam_1t^{k_1},...,
\lam_nt^{k_n})=t^l F^\del(\lam_1,...,\lam_n)=0$,
$l=i_1k_1+...+i_nk_n$, $(i_1,...,i_n)\in\del$,
and there exist smooth functions $\varphi_1(t)=\lam_1+O(t)$,
..., $\varphi_n=\lam_n+O(t)$ such that the curve
$\gam(t)=(t^{k_1}\varphi_1(t),...,t^{k_n}\varphi_n(t))$ lies on
$\{F=0\}$. Indeed,
$$F(\gam(t))=t^\lam F^\del(\varphi_1,...,\varphi_n)+
t^{\lam+1}G(t,\varphi_1,...,\varphi_2)=0$$
$$\Longleftrightarrow\quad F^\del(\varphi_1,...,\varphi_n)+
t G(t,\varphi_1,...,\varphi_2)=0\ ,$$
hence the existence of $\varphi_1,...,\varphi_n$ defined for
small $t>0$ follows from the implicit function theorem and
the fact that $F^\del=0$ is nonsingular in $(\C^*)^n$. The above
argument shows that
$$\tmu_\del(\lam_1,...,\lam_n)=\lim_{t\to 0}\tmu_\Del(\gam(t))\ .$$

(iv) Assume that $\dim\Del=n$.
Let $\del\subset\Del$ be a proper face, $\dim\del=s>0$, and
$\ow=(y_1v^0_1,...,y_nv^0_n)\in\C I(\del)\cap\C Ch(F)$, where
$(y_1,...,y_n)\in I(\del)$, $(v^0_1,...,v^0_n)\in(S^1)^n$.
We will describe $\C Ch(F)$ in a neighborhood~$U_\ow$ of
the point~$\ow$ in~$\C\Del$.

Represent $\R^n$ as $\R^s\times\R^{n-s}$ and denote by
$\pr_s:\R^n\to\R^s$, $\pr_{n-s}:\R^n\to\R^{n-s}$ the natural
projections.
According to Proposition \ref{p7}, we can assume that $\pr_{n-s}
(\del)$ is a point $(p_{s+1},...,p_n)$ with positive integers
$p_{s+1},...,p_n$, so $\widetilde\del=\pr_s(\del)=\del
-(0,...,0,p_{s+1},...,p_n)$, and
$\widetilde\Del\backslash\widetilde\del$
lies in the positive orthant.
Then
$F(z_1,...,z_n)=
z_{s+1}^{p_{s+1}}...z_n^{p_n}\widetilde F(z_1,...,z_n)$, where
\begin{equation}
\widetilde F(z_1,...,z_n)=\widetilde F^{\widetilde\del}
(z_1,...,z_s)+
\sum_{(i_1,...,i_n)\in\Lam}
A_{i_1...i_n}z_1^{i_1}...z_s^{i_s}\cdot M_{i_{s+1}...i_n}\ ,
\label{e16}\end{equation}
$$\Lam=\{(i_1,...i_n)\in\widetilde\Del\ :
\ i_{s+1},...,i_n>0\}\cap\Z^n,
\quad M_{i_{s+1}...i_n}=z_{s+1}^{i_{s+1}}...z_n^{i_n}\ .$$

Note that there are uniquely defined $x^0_1,...,x^0_s>0$ such that
$(y_1,...,y_s)=\mu_{\widetilde\del}(x^0_1,...,x^0_s)$.
The condition on $\oz=(\ox,\ov)\in(\R^*_+)^n\times(S^1)^n=(\C^*)^n$
such that $\tmu_\Del(\ox,\ov)\in U_\ow$ can be expressed as
follows. Clearly, for $\ov=(v_1,...,v_n)$, it is
\begin{equation}
|v_1-v^0_1|<\eps,\ ...,\ |v_n-v^0_n|<\eps,\label{e14}
\end{equation}
where $\eps>0$ is sufficiently small. For $\ox=(x_1,...,x_n)$, we have
that $\mu_{\widetilde\Del}(\ox)$, equal to
$$\frac{\sum_{(i_1,...,i_s)\in\widetilde\del}x_1^{i_1}...
x_s^{i_s}\cdot(i_1,...,i_s,0,...,0)+
\sum_{(i_1,...,i_n\in\Lam}x_1^{i_1}...
x_s^{i_s}\cdot|M_{i_{s+1}...i_n}|\cdot(i_1,...,i_n)}
{\sum_{(i_1,...,i_s)\in\widetilde\del}x_1^{i_1}...
x_s^{i_s}+
\sum_{(i_1,...,i_n\in\Lam}x_1^{i_1}...
x_s^{i_s}\cdot|M_{i_{s+1}...i_n}|}\ ,$$
is close to $(y_1,...,y_s,0,...,0)$.
This means that in the latest
expression the summands over $\Lam$ are small with respect to
the sums over $\widetilde\del$; hence $\pr_s(\mu_{\widetilde\Del}
(\ox))$ is close to $\mu_{\widetilde\del}(x_1,...,x_s)$, which
means
that $(x_1,...,x_s)$ is close to $(x^0_1,...,x^0_s)$. Coming back
to the sums over $\Lam$, we obtain that $M_{i_{s+1}...i_n}$
must be close to zero, $(i_{s+1},...,i_n)\in\pr_{n-s}(\Lam)$.
So, without loss of generality the conditions on $\ox$ can be
expressed as
\begin{equation}
|x_l-x^0_l|<\eps,\ l=1,...,s,\quad
x_{s+1}^{i_{s+1}}...x_n^{i_n}<\eps,\
(i_{s+1},...,i_n)\in\pr_{n-s}(\Lam)\ .\label{e15}
\end{equation}
Let $U'\subset(\R^*_+)^n$ and $U''\subset(S^1)^n$ be given
by~(\ref{e15}) and~(\ref{e14}), respectively,
and note that $U'=\pr_s(U')\times\pr_{n-s}(U')$.
Put $\sig=\pr_{n-s}(\widetilde\Del)$. Then $\mu_\sig$ takes
$\pr_{n-s}(U')$ onto a set $U_0\backslash\partial\sig$, where
$U_0$ is a neighborhood of the origin (a vertex of $\sig$~!) in
$\sig$.

Since $\widetilde F$ is non-degenerate, we can assume that
$$\frac{\partial\widetilde F}{\partial x_1}(x^0_1,...,x^0_s,
0,...,0)=\frac{\partial\widetilde F^{\widetilde\del}} {\partial
x_1}(x^0_1,...,x^0_s) \ne0\ .$$ Hence, by (\ref{e16}) and the
implicit function theorem, the hypersurface $\{F=0\}$ can be
described in $U'\times U''$ (with, possibly, smaller $\eps$) by
equations
\begin{equation}
x_1=x^0_1+G,\quad v_1=v^0_1+H\ ,\label{e17}
\end{equation}
where~$G$ and~$H$ are
vanishing at zero smooth functions of $x_2-x^0_2$, ...,
$x_s-x^0_s$, $v_2-v^0_2$, ...,
$v_n-v^0_n$, and $M_{i_{s+1}...i_n}$,
$(i_{s+1},...,i_n)\in\pr_{n-s}(\Lam)$,
where $x_2,...,x_n$, $v_2,...,v_n$ are arbitrary
satisfying~(\ref{e14}) and~(\ref{e15}).

Now we show that $U_\ow$ is homeomorphic to
$\pr_s(U')\times U_0\times U''$ via an extension of $\tmu_\Del$.
Indeed,
$U_\ow\backslash\partial\C\Del$ is homeomorphic via $\tmu_\Del$ to
$\pr_s(U')\times\pr_{n-s}(U')\times U''\simeq
\pr_s(U'')\times(U_0\backslash\partial\sig)\times U''$.
As we saw, the points
in $U_\ow\cap\partial\Del$ are $\lim_{t\to 0}\tmu_\Del(\ox(t),
\ov(t))$
for all curves $\ox(t)\in U'$, $\ov(t)\in U''$, $t>0$,
such that $\lim_{t\to 0}\ov(t)\in U''$,
$\lim_{t\to 0}(x_1(t),...,x_s(t))\in\pr_s(U')$ and
$$x_{s+1}(t)=\lam_1t^{k_1}+O(t^{k_1+1}),\ ...,\ x_n(t)=
\lam_{n-s}t^{k_{n-s}}+O(t^{k_{n-s}+1})\ .$$
Restriction (\ref{e15}) means,
first, that the vector $(k_1,...,k_{n-s})$
belongs to the cone generated by the interior normal vectors
to the $(s-1)$-dimensional
faces of $\sig$, which contain the origin, and, second, that
$\lam_1^{i_{s+1}}...\lam_{n-s}^{i_n}<\eps$ for all
$(i_{s+1},...,i_n)\in\pr_{n-s}(\Lam)$ orthogonal
to $(k_1,...,k_{n-s})$.
On the other hand, the same data determine a point
$\lim_{t\to 0}\mu_\sig(x_{s+1}(t),...,x_n(t))\in U_0\cap
\partial\sig$. This gives a one-to-one correspondence between
$U_\ow\cap\partial\Del$ and $\pr_s(U')\times(U_0\cap
\partial\sig)\times U''$,
which respects the combinatorial structure: the
interior of a face $\psi\subset U_\ow\cap\partial\Del$
corresponds to $\pr_s(U')\times\pr_{n-s}(I\widetilde\psi)
\times U''$, where $\widetilde\psi=\psi-(0,...,0,p_{s+1},
...,p_n)$. The above correspondence
provides a homeomorphism of
$U_\ow$ and $\pr_s(U')\times U_0\times U''$, since close points
in $U_\ow\cap\partial\Del$ come from close curves
$\gam(t)$ which in turn bear close points in
$\pr_s(U')\times(U_0\cap
\partial\sig)\times U''$. Consequently, we can introduce in
$U_0$ coordinates $\tet_1,...,\tet_{n-s}$ such that
$$U_0=\{0\le\tet_1<\eps,\ -\eps<\tet_j<\eps,\ j=2,...,n-s\},
\quad U_0\cap\partial\sig=\{\tet_1=0\}\ ,$$
and $M_{i_{s+1}...i_n}$, $(i_{s+1},...,i_n)\in\pr_{n-s}(\Lam)$,
are continuous functions of $\tet_1,...,\tet_{n-s}$,
vanishing at zero. Thus, formulae (\ref{e17})
allow us to parameterize
$\C Ch(F)\cap U_\ow$ by
the following subset in $\pr_s(U')\times U_0\times U''$:
$$x_1=x^0_1+G',\quad v_1=v^0_1+H'\ ,$$
where~$G'$ and~$H'$ are vanishing at zero smooth functions of
$x_2,...,x_s,v_2,...,v_n,\tet_1,...,\tet_{n-s}$.
This subset is a subvariety of
codimension 2 with a boundary on
$\pr_s(U')\times(U_0\cap\partial\sig)
\times U''$.
Hence $\C Ch(F)\cap U_\ow$ is a subvariety of codimension 2 in $U_\ow$
with boundary on $U_\ow\cap\partial\Del$.

(v) The statement of Proposition for $\R Ch(F)$ can be proved
as in the complex case.
\proofend

\begin{corollary}\label{c1}
If $\dim\Del=2$, $\Del\subset\R^2$, then $\C Ch(F)$ is a smooth
surface with boundary.
\end{corollary}

\subsection{Digression: real resolution of toric varieties}
Let $\Del\subset(\R^*_+)^n$ be an $n$-dimensional polytope, $F$ be
a non-degenerate polynomial with Newton polytope $\Del$. Denote by
$\Tor_\C(\Del)$ the toric variety over $\C$ defined by $\Del$, and
by $Z(F)$ the (closed) hypersurface in $\Tor_\C(\Del)$ defined by
$F$.

\begin{proposition}\label{p14}
There exists an
equivariant surjective map $\nu:\C\Del\to \Tor_\C(\Del)$ such that
\begin{itemize}
\item $\nu(\C Ch(F))=Z(F)$,
\item $\nu\big|_{\C I(\Del)}$ takes $\C I(\Del)$ diffeomorphically
to $(\C^*)^n\subset \Tor_\C(\Del)$,
\item for any face $\del$, $\dim\del=s<n$, the map
$\nu\big|_{\C I(\del)}$ is a submersion of $\C I(\del)$ onto
$(\C^*)^s\subset \Tor_\C(\del)\subset \Tor_\C(\Del)$.
\end{itemize}
\end{proposition}

{\bf Proof.}
Let us represent $\Tor_\C(\Del)$ as the closure of
the variety
$$\{(z_1^{i_1}...z_n^{i_n})_{(i_1,...,i_n)\in\Del\cap\Z^n}
\ :\ (z_1,...,z_n)\in(\C^*)^n\}
\subset\C P^{N-1}\ ,$$
where $N=\#(\Del\cap\Z^n)$ (see, for example \cite{Fu}).
We define $\nu\big|_{\C I(\Del)}=(\tmu_\Del)^{-1}$, which is a
required diffeomorphism of $\C I(\Del)$ and
$(\C^*)^n \subset \Tor_\C(\Del)$, and extend it to $\partial\C\Del$.
Namely,
given a face $\del\subset\Del$ and a point $\ow=\tmu_\del(\oz)\in
\C I(\del)$, $\oz=(z_1,...,z_n)\in(\C^*)^n$, we put
$\nu(\ow)=(a_{i_1...i_n})_{(i_1,...,i_n)\in\Del\cap\Z^n}$, where
$$a_{i_1...i_n}=z_1^{i_1}...z_n^{i_n},\ (i_1,...,i_n)\in\del,
\quad a_{i_1...i_n}=0,\ (i_1,...,i_n)\in\Del\backslash\del\ .$$
Clearly, $\nu\big|_{\C I(\del)}$ satisfies the required property.
It remains only to explain that $\nu$ is continuous. Indeed,
in the previous notation
$\ow=\lim_{t\to\infty}\tmu_\Del(\gam(t))$, where
$$\gam(t)=(z_1t^{k_1}+O(t^{k_1-1}),...,z_nt^{k_n}+O(t^{k_n-1})),
\quad t>0\ ,$$
is a curve in $(\C^*)^n$ with $(k_1,...,k_n)\in\Z^n\backslash\{0\}$
belonging to the dual cone of $\del$ with respect to $\Del$.
Then one can easily check that $\lim_{t\to\infty}\gam(t)\in
\C P^{N-1}$ is just $\nu(\ow)$ defined as above.
\proofend

\subsection{Gluing of charts}\label{s2}
Let ${\cal S}$ be a subdivision of a polytope $\Del\subset\R^n$
into equidimensional polytopes: $\Del=\Del_1\cup...\cup\Del_N$
({\it i.e.}, $\Del_i\cap\Del_j$ is empty or a common proper face),
and let ${\cal A}=\{A_\oi,\ \oi \in\Del\cap\Z^n\}$ be a collection
of complex numbers such that $A_\oi\ne 0$ if $\oi$ is a vertex of
some $\Del_i$, $1\le i\le N$. Further on, speaking on subdivisions
of polytopes and corresponding collections of numbers, we always
assume the above properties.

Assume that the polynomials
$$F_k(z_1,...,z_n)=\sum_{(i_1,...,i_n)\in\Del_k}
A_{i_1...i_n}z_1^{i_1}\cdot...\cdot z_n^{i_n}\ ,\quad k=1,...,N,$$
are non-degenerate. The union of the complex charts of the
polynomials $F_1,...,F_N$ $$\C Ch({\cal S},{\cal
A})=\bigcup_{k=1}^N\C Ch(F_k)$$ is called a {\it C-hypersurface in
$\C\Del$}. The union of the real charts of $F_1,...,F_N$ $$\R
Ch({\cal S},{\cal A})=\bigcup_{k=1}^N\R Ch(F_k)$$ is called a {\it
C-hypersurface in $\R\Del$}.

\begin{proposition}\label{l9} Let $\Del\subset(\R^*_+)^n$.
Then the set $\C Ch({\cal S},{\cal A})$ is a PL-manifold
of codimension 2 in $\C\Del$ with boundary
$\partial\C Ch({\cal S},{\cal A})=
\partial\C\Del\cap\C Ch({\cal S},{\cal A})$.
If all the numbers in~${\cal A}$ are real,
then $\C Ch({\cal S},{\cal A})$ is
invariant with respect to $\conj$, the set
$\R Ch({\cal S},{\cal A})$ is a PL-manifold
of codimension 1 in $\R\Del$ with boundary
$\partial\R Ch({\cal S},{\cal A})=
\partial\R\Del\cap\R Ch({\cal S},{\cal A})$.
\end{proposition}

\begin{remark}\label{r2}

(a) Neither $\C Ch({\cal S},{\cal A})$, nor $\R Ch({\cal S},{\cal A})$ are
smooth in general. Indeed, for the polynomials
$$F_1(x,y)=x(y-x-1),\quad F_2(x,y)=xy-x+1$$
with the Newton
triangles
$$\Del_1=[(1,0),(1,1),(2,0)],\quad\Del_2=[(0,0),(1,0),(1,1)]$$
the closures of curves
$$\mu_{\Del_1}(\{F_1=0\})=\left\{\left(\frac{2+3x}{2+2x}\ ,
\ \frac{x+1}{2+2x}\right)\ :\ x>0\right\},$$
$$\mu_{\Del_2}(\{F_2=0\})=\left\{\left(\frac{2x-1}{2x}\ ,
\ \frac{x-1}{2x}\right)\ :\ x>0\right\}$$
have the common limit point $(1,1/2)$ on the edge $[(1,0),(1,1)]$,
and different tangent lines
$$y=\frac12,\quad y=x-\frac12$$
at this point.

(b) Similarly to Proposition
\ref{l1}, the statement of Proposition \ref{l9}
holds for $\Del$ intersecting coordinate hyperplanes,
when substituting
$\C\Del\backslash\Sing(\C\Del)$ and $\R\Del\backslash\Sing(\C\Del)$
for $\C\Del$ and $\R\Del$, respectively.
\end{remark}

{\bf Proof of Proposition \ref{l9}.}
In view of Proposition \ref{l1}, we have to study only the behavior of
$\C Ch({\cal S},{\cal A})$ of faces $\del\subset\Del_k\cap\Del_l$,
$k\ne l$, $\dim\del=s>0$.

Let $I(\del)\subset I(\Del)$,
$\del=\Del_1\cap...\cap\Del_k$, $\del\not\subset\Del_{k+1}\cup...
\cup\Del_N$. Pick a point $w\in\C I(\del)\cap
\C Ch({\cal S},{\cal A})$. Following the proof of Proposition
\ref{l1},
we impose assumptions of step
(iv) and obtain that a neighborhood $U_w$ of $w$ in
$\C\Del$ is the union of neighborhoods $U_{w,1},...,U_{w,k}$
of $w$ in $\C\Del_1,...,\C\Del_k$, respectively, so that
$U_{w,m}$ is parameterized by $\pr_s(U'_m)\times U_{0,m}
\times U''$, where $\pr_s(U'_1)=...=\pr_s(U'_k)$,
$U''$ is common for all $m=1,...,k$, and $U_{0,m}$ is a neighborhood
of the point $\pr_{n-s}(\del)=(p_{s+1},...,p_n)$ in
$\pr_{n-s}(\Del_m)$, and, moreover, these parameterizations are
compatible on the common faces of $\Del_1,...,\Del_k$. In turn,
$\C Ch({\cal S},{\cal A})\cap U_{w,m}$ is parameterized by
$$\prod_{j=2}^s\{|x_j-x^0_j|<\eps\}\times U_{0,m}\times
\prod_{j=2}^n\{|v_j-v^0_j|<\eps\}\ ,$$
and by (\ref{e13}) these parameterizations are compatible along
common faces of $\Del_1,...,\Del_k$. Hence
$\C Ch({\cal S},{\cal A})\cap U_w=\bigcup_m\C Ch({\cal S},
{\cal A})\cap U_{0,m}$ is parameterized by
$$\prod_{j=2}^s\{|x_j-x^0_j|
<\eps\}\times\bigcup_{m=1}^k
U_{0,m}\times
\prod_{j=2}^n\{|v_j-v^0_j|<\eps\}\ ,$$
which is homeomorphic to a $(2n-2)$-ball.

Similarly one treats the case $\del\subset\partial\Del$ and
the real case.
\proofend

\begin{definition}\label{d4}
A homeomorphism of (resp., an isotopy in)
$\C\Del$
is called {\it tame} if
for any face~$\del$
of $\Del$
the restriction of this homeomorphism
(resp., isotopy) to $\C\del$
is a homeomorphism of (resp., an isotopy in) $\C\del$.
In addition, we call such objects
{\it equivariant} if they commute with $\conj$.
\end{definition}

\begin{remark}\label{r4}
For given $\Del$ and $\cal S$, and different ${\cal A}, {\cal
A}':\Del\cap\Z^n \to \C$, the C-hypersurfaces $\C Ch({\cal
S},{\cal A})$ and $\C Ch({\cal S},{\cal A}')$ are tame isotopic in
$\C\Del$ (not equivariantly, in general). Indeed, one can connect
${\cal A}$ and ${\cal A}'$ by a family ${\cal A}_t$, $t\in[0,1]$
such that the polynomials $F_{k,t}$ are non-degenerate for all
$k=1,...,N$, $t\in[0,1]$.
\end{remark}

\subsection{Projectivization of charts}\label{s1}
Let~$d$ and~$n$ be positive integers. Denote by $T_d^n$
the simplex in
$\R^n$ with vertices $(0,0,...,0)$, $(d,0,...,0)$, ...,
$(0,...,0,d)$. Observe that $\C T_d^n$ is homeomorphic to a closed
ball in $\C^n$. The group $S^1=\{|v|=1\}$ acts freely on $\partial\C
T_d^n$ by
$$v\in S^1,\ (z_1,...,z_n)\in\partial
\C T_d^n\ \mapsto\ (z_1v,...,z_nv)
\in\partial\C T_d^n\ .$$

\begin{proposition}\label{p1}
The quotient $\C T_d^n/S^1$ is
equivariantly homeomorphic to
the projective space
$\C P^n$ so that the faces of $T_d^n$ naturally correspond to
the coordinate planes in $\C P^n$.
\end{proposition}

{\bf Proof.}
The extended moment map $\tmu_{T_d^n}$ takes $\C^n$ diffeomorphically
onto $\Int(\C T^n_d)$. Its inverse
extends up to a surjective map
$\nu_d^n:\C T_d^n\to\C P^n$ in the following way.
Consider a point
$(0,z_1,...,z_n)\in\C P^n\backslash\C^n$, the line
$\gam(t)=(z_1t,...,z_nt)\in\C^n$, $t\in(0,\infty)$,
and $\lim_{t\to\infty}\tmu_{T_d^n}(\gam(t))$.
Then
$\lim_{t\to\infty}\tmu_{T_d^n}(\gam(t))=\tmu_\del(z_1,...,z_n)$,
where
$$\del=\left\{\sum_{z_k\ne 0}i_k=d\right\}\cap\bigcap_{z_k=0}
\{i_k=0\}$$
is a face of $T_d^n$. Then we put $\nu_d^n(\tmu_\del(z_1,...,z_n)=
(0,z_1,...,z_n)\in\C P^n$. One can easily verify that
$\nu_d^n$ is well-defined, surjective, continuous, commutes
with $\conj$, sends orbits of the $S^1$-action
on $\partial\C T_d^n$ into the points of the hyperplane $\{z_0=0\}$
in $\C P^n$, and establishes a one-to-one correspondence between
the faces of $T^n_d$ and the coordinate planes in $\C P^n$.
\proofend

\begin{remark}\label{r5}
Further on we always assume $\C T_d^n/S^1$ to be identified with
$\C P^n$ via the map $\nu_d^n$.
\end{remark}

\begin{definition}\label{d1}
In the notation of section \ref{s2}, given a subdivision
$T_d^n=\Del_1\cup...\cup\Del_N$ and a collection of complex numbers
${\cal A}=\{A_\oi\ :\ \oi \in T_d^n\cap\Z^n\}$, define
a {\it C-hypersurface}
{\it of degree}~$d$
in~$\C P^n$ as
$$P\C Ch({\cal S},{\cal A})=\nu_d^n(\C Ch({\cal S},{\cal A}))
\subset\C P^n\ .$$
If the numbers ${\cal A}=\{A_\oi\ :\ \oi \in
T_d^n\cap\Z^n\}$ are real, the C-hypersurface is called {\it real}
and its real part $\nu_d^n(\R Ch({\cal S},{\cal A}))$
is denoted by $P\R Ch({\cal S},{\cal A})$.
\end{definition}

\begin{proposition}\label{p8}
If $F$ is a non-degenerate polynomial of degree $d$ in $n$
variables, then $P\C Ch(F)\subset\C P^n$ coincides with the
projective closure $P\{F=0\}$ of the affine hypersurface
$\{F=0\}$.
\end{proposition}

{\bf Proof.} Straightforward from $\nu_d^n(\C Ch(F))
=P\{F=0\}$ and the fact that $\C Ch(F)\cap\partial(\C T_d^n)$
is invariant with respect to the $S^1$-action.
\proofend

\begin{definition}\label{d2}
A homeomorphism of (resp., an isotopy in)
$\C P^n$
is called {\it tame} if
its restriction to any coordinate plane
is a homeomorphism of (resp., an isotopy in)
that plane.
In addition, we call such objects
{\it equivariant} if they commute with $\conj$.
\end{definition}

\begin{proposition}\label{l3}
In the previous notation,
$P\C Ch({\cal S},{\cal A})$ is a PL-submanifold
of codimension~$2$ in $\C P^n$.
It is invariant with respect to $\conj$ if
all the numbers in~${\cal A}$ are real.
\end{proposition}

{\bf Proof.}
First, we show that
the C-hypersurface
$\C Ch({\cal S},{\cal A})$ is a PL-manifold in $\C T_d^n$
of codimension 2 with boundary
$$\partial\C Ch({\cal S},{\cal A})=\C Ch({\cal S},{\cal A})\cap
\partial\C T_d^n\ ,$$
and if all the numbers in~${\cal A}$ are real
$\C Ch({\cal S}, {\cal A}$
is invariant with respect to
$\conj$, and $\Fix(\conj)=\R Ch({\cal S},{\cal A})$ is a
PL-manifold in $\R T_d^n$
of codimension 1 with boundary
$$\partial\R Ch({\cal S},{\cal A})=\R Ch({\cal S},{\cal A})\cap
\partial\R T_d^n\ .$$

By Proposition \ref{l9} it is enough to consider
only the union of charts along the coordinate hyperplanes.
For any polytope $\Del\subset\R^n$,
put $\Del_t=(t,...,t)+\Del$, $t\ge 0$.

The family $\tmu_{\Del_{k,t}}$, $0\le t\le 1$, connects
$\tmu_{\Del_k}$ with $\tmu_{\Del_{k,1}}$, $k=1,...,N$, and defines
an isotopy of $\Cl(\tmu_{\Del_{k,t}}(F_k=0))$,
the closure of $\tmu_{\Del_{k,t}}(F_k=0)$, $0<t\le 1$,
which degenerates into $\C Ch(F_k)$, $k=1,...,N$.
The same holds for any proper face $\del$ of $\Del_1,...,
\Del_N$ and the corresponding truncation of $F_1,...,F_N$.

Since $\Sing(\C\Del_{k,1})=\emptyset$,
$k=1,...,N$, the set
${\cal M}_1=\bigcup_{k=1}^N\C Ch(z_1...z_nF_k)$
is a PL-manifold with boundary on $\partial\C \widetilde T_d^n$,
according to Proposition \ref{l1}
(where $\widetilde T_d^n=(1,...,1)+T_d^n$).
Denote some faces of $\widetilde T_d^n$ as follows:
$$\widetilde T_d^n(k)=\widetilde T_d^n\cap\{i_k=1\},\quad
\widetilde T_d^n(k_1,...,k_s)=\bigcap_{j=1}^s
\widetilde T_d^n(k_j)\ .$$

Put $\del=\widetilde T_d^n(k)$. By (\ref{e12})
\begin{equation}\C\del
\cong\C T_d^{n-1}\times S^1.\label{e10}
\end{equation}
Similarly, by (\ref{e12})
${\cal M}_1\cap\C\del$ is a
PL-manifold of dimension $2n-3$ with boundary
on $\partial\C\del$, satisfying
\begin{equation}
{\cal M}_1\cap
\C\del\cong\left(\bigcup_{l=1}^N\C Ch(z_1...z_nF_l^{\del})
\right)\times S^1,\label{e7}
\end{equation}
where the product structure is compatible with that in (\ref{e10}).
Using
the product representations (\ref{e10}) and (\ref{e7}), we define
$$D^n(k)=\{(w_1,...,w_n)\ :\ (w_1,...,w_{k-1},w_{k+1},...,w_n)\in
\C T_d^{n-1},\ |w_k|\le 1\}$$
$$\cong\C T_d^{n-1}\times D^2\ ,$$
where $D^2$ is
the closed 2-dimensional unit ball, and inside $D^n(k)$
the $(2n-2)$-dimensional manifold $M(k)$:
$$\{(w_1,...,w_n)\ :\ (w_1,...,w_{k-1},1,w_{k+1},...,w_n)
\in\bigcup_{l=1}^N\C Ch(z_1...z_nF_l^{\del}),\ |w_k|\le 1\}$$
$$\cong\left(\bigcup_{l=1}^N\C Ch(\widetilde F_l^{\del})
\right)\times D^2$$
with boundary on $\partial D^n(k)$ and such that
$M(k)\cap\C\del=\bigcup_{l=1}^N\C Ch(z_1...z_nF_l^{\del})$.

Similarly, given $\del=\widetilde T_d^n
(k_1,...,k_s)$, we have by (\ref{e12}),
$$\C\del\cong\C T_d^{n-s}\times(S^1)^s,\quad
{\cal M}_1\cap\C\del
\cong\left(\bigcup_{l=1}^N\C Ch(z_1...z_nF_l^{\del})
\right)\times(S^1)^s,$$
and one defines
$$D^n(k_1,...,k_s)=
\{(w_1,...,w_n)\ :\ (w_j)_{j\ne k_1,...,k_s}\in
\C T_d^{n-s},\ |w_{k_1}|,...,|w_{k_s}|\le 1\}$$
$$\cong\C T_d^{n-s}\times(D^2)^s\ ,$$
and the $(2n-2)$-dimensional manifold $M(k_1,...,k_s)$:
$$\{(w_1,...,w_n)\ :\ (w_j)_{j\ne k_1,...,k_s}
\in\bigcup_{l=1}^N\C Ch(z_1...z_nF_l^{\del}),\
|w_{k_1}|,...,|w_{k_s}|\le 1\}$$
$$\cong\left(\bigcup_{l=1}^N\C Ch(z_1...z_nF_l^{\del})
\right)\times(D^2)^s$$
with boundary on $\partial D^n(k_1,...,k_s)$ and such that
$$M(k_1,...,k_s)\cap\C\del=
\bigcup_{l=1}^N\C Ch(z_1...z_nF_l^{\del})\ .$$

It can easily be shown that
${\cal T}_1=\C\widetilde T_d^n
\cup\bigcup_{s=1}^n\bigcup_{1\le k_1<...<
k_s}D^n(k_1,...,k_s)$
is the convex hull of $\C\widetilde T_d^n$ in $\C^n$, and
${\cal M}_1\cup\bigcup_{s=1}^n\bigcup_{1\le k_1<...<
k_s}M(k_1,...,k_s)$
is a PL-manifold in ${\cal T}_1$ of codimension 2 with boundary
on $\partial{\cal T}_1$.

Now, for
any $t\in(0,1)$, in the same way one constructs similar objects:
\begin{itemize}
\item complexifications of polytope
$$\C\Del_{k,t}=\Cl(\tmu_{\Del_k,t}((\C^*)^n)),\quad\C T_d^n(t)=
\bigcup_{k=1}^N\C\Del_{k,t},$$
\item union of charts
${\cal M}_t=\bigcup_{k=1}^N\Cl(\tmu_{\Del_k,t}(\{F_k=0\}))$,
which is
a PL-manifold of codimension 2 in $\C T_d^n(t)$
with boundary on $\partial\C
T_d^n(t)$,
\item completion of $\C T_d^n(t)$ up to its convex hull
$${\cal T}_t=\C T_d^n(t)
\cup\bigcup_{s=1}^n\bigcup_{1\le k_1<...<
k_s}D_t^n(k_1,...,k_s),$$
where
$$D_t^n(k_1,...,k_s)=
\{(w_1,...,w_n)\ :\ (w_j)_{j\ne k_1,...,k_s}\in
\C T_d^{n-s},\ |w_{k_1}|,...,|w_{k_s}|\le t\}$$
$$\cong\C T_d^{n-s}\times(D^2(t))^s\ ,$$
$D^2(t)$ is a disc of radius $t$,
\item
a $(2n-2)$-dimensional manifold
$$M_t(k_1,...,k_s)=\Big\{(w_1,...,w_n)\ :$$
$$\ (w_j)_{j\ne k_1,...,k_s}
\in\bigcup_{l=1}^N\Cl(\tmu_{\del_t}(\{F_l^{\del}=0\})),\
|w_{k_1}|,...,|w_{k_s}|\le t\Big\}$$
$$\cong\left(\bigcup_{l=1}^N\Cl(\tmu_{\del_t}(\{F_l^{\del}
=0\}))
\right)\times(D^2(t))^s$$
\item a PL-manifold
\begin{equation}
{\cal M}_t\cup\bigcup_{s=1}^n\bigcup_{1\le k_1<...<
k_s}M_t(k_1,...,k_s)\label{e8}
\end{equation}
of codimension 2
in ${\cal T}_t$ with boundary on $\partial{\cal T}_t$.
\end{itemize}

If $t$ varies from 1 to 0, ${\cal T}_t$ contracts from ${\cal T}_1$
to $\C T_d^n$ and the manifold (\ref{e8}) naturally contracts
into the C-hypersurface $\C Ch({\cal S},{\cal A})$ with boundary
on $\partial\C T^n$.

To complete the proof of Proposition \ref{l3}
we have to verify only that $\partial\C Ch({\cal S},{\cal A})$
is invariant with respect to the $S^1$-action on
$\partial\C T_d^n$, and
the quotient $\C Ch({\cal S},{\cal A})/S^1$ is a
closed manifold.

First, note
that $\partial\C Ch({\cal S},{\cal A})$
is invariant with respect to the $S^1$-action. Indeed,
$\partial\C Ch({\cal S},{\cal A})$ consists of charts of
homogeneous polynomials which are invariant with
respect to the $S^1$-action on
$\partial\C T_d^n$. Second,
for any edge of $T_d^n$
there exists a combination of an automorphism of $\Z^n$ and shifts
which puts this edge on a coordinate axis and the adjacent faces
of $T_d^n$ on the corresponding coordinate planes. Since
such a transformation is compatible with $\nu_d^n$ defined
via the extended moment map, the question on the behavior of
$\C Ch({\cal S},{\cal A})$ on $\partial\C T_d^n$ is
reduced to that on coordinate planes, which has been treated above.
\proofend

\begin{remark}
Similarly to Remark \ref{r4},
two C-hypersurfaces in $\C P^n$ of the same degree and
with the same subdivision ${\cal S}$ of $T_d^n$ are tame
isotopic in $\C P^n$ (not equivariantly, in general).
\end{remark}

\section{Topology of complex C-hypersurfaces}\label{compl-hyper}

\subsection{Complex version of the Viro theorem}\label{Viro}

\begin{theorem}\label{l5}
{\rm (Affine complex Viro theorem, see~\cite{Vi1})}.
In the notation of section \ref{s2}, if $\Del\subset(\R^*_+)^n$
and a subdivision ${\cal S}\ : \ \Del=\Del_1\cup...\cup\Del_N$ is
defined by a convex piecewise-linear function $\nu:T_d^n\to\R$,
then $\C Ch({\cal S},{\cal A})$ is tame isotopic in $\C\Del$
(equivariantly, if all the numbers in~${\cal A}$ are real) to $\C
Ch(F)$, where $$F(z_1,...,z_n,t)=
\sum_{(i_1,...,i_n)\in\Del}A_{i_1...i_n}t^{\nu(i_1,...,i_n)}
z_1^{i_1}\cdot... \cdot z_n^{i_n}$$ and $t=\const>0$ is
sufficiently small.
\end{theorem}

{\bf Proof.}
Without loss of generality one can suppose that
$\dim\Del=n$, and $\nu$ is positive
integral-valued at
integral points. Consider the polynomial
$\hat F(z_1,...,z_n,t)=(1+t)F(z_1,...,z_n,t)$ as a
polynomial in
$t,x_1,...,x_n$, and denote by $\hat\Del$ its Newton polytope
in $\R^{n+1}$. The graph of $\nu$ is the ``lower" part of
$\partial\hat\Del$, naturally projected onto $\Del\subset\R^n=
\{i_{n+1}=0\}\subset\R^{n+1}$. Denote by $\hat\Del_k$
the part of the graph of $\nu$ projected onto
$\Del_k$, $k=1,...,N$. Denote by $B$ the part of
$\hat\Del$ projected onto $\partial\Del$. Clearly,
$B$ can be viewed as $\partial\Del\times[0,1]$.

{\it Step 1}. Introduce the set
$$H_0=\{(w_1,...,w_n,\nu(|w_1|,...,|w_n|))\ :\ (|w_1|,...,
|w_n|)\in\Del\}=\bigcup_{k=1}^N\hat\Del_k
\times(S^1)^n\subset\partial\C\hat\Del\ ,$$
and the family of hypersurfaces
$H_a=\Cl(\tmu_{\hat\Del}(\{t=a\}))\subset\C\hat\Del$, $a\in(0,\eps)$.
We will show that
the family $H_a$, $a\in[0,\eps)$, is an isotopy in $\C\hat\Del$.

The map $\tmu_{\hat\Del}$ takes $\{t=a\}\subset(\C^*)^{n+1}$
diffeomorphically into $\C I(\hat\Del)$.
Considering $\lim_{\tau\to 0}\mu_{\hat\Del}(\gam(\tau))$
for all curves
$$\gam(\tau)=\left(\lam_1x_1^{k_1},...,\lam_nx_n^{k_n},a\right)
\subset(\R^*_+)^{n+1},
\quad\tau>0,$$
with $\lam_1,...,\lam_n>0$, $(k_1,...,k_n)\ne 0$, one can show
(as in the proof of Proposition \ref{l1}) that $H_a$ being the closure
of $\tmu_{\hat\Del}(\{t=a\})$ is a PL-manifold with boundary
$\partial H_a=H_a\cap B$, which projects onto $\partial\Del$.
Moreover, $H_a\cap H_b=\emptyset$, $a\ne b$.
The projection on the first $n$ coordinates
provides a homeomorphism of $H_a$, $a>0$, and
$\C\Del$. Indeed, in $\C I(\hat\Del)$
this fact is reduced to the claim that
$$(x_1,...,x_n)\in(\R^*_+)^n\mapsto
\frac{\sum_{(i_1,...,i_n)\in\Del}x_1^{i_1}...x_n^{i_n}a^{\nu(i_1,...,
i_n)}\cdot(i_1,...,i_n)}{\sum_{(i_1,...,i_n)\in\Del}
x_1^{i_1}...x_n^{i_n}a^{\nu(i_1,..., i_n)}}\in I(\Del)$$
is a homeomorphism, which holds since the above map is
a kind of the moment map
\cite{A1, A2}. Similarly one shows that the projection is one-to-one
on $\partial H_a\subset B$, taking into account that the limit of
$\mu_{\hat\Del}$ on a face of $B$ is the moment map of this face.
To finish the proof one should note that
$\bigcup_{a>0}H_a$ fill the intersection of $\C I(\hat\Del)$ with
the space $\im(w_{n+1})=0$, where $w_1,...,w_{n+1}$ are
coordinates in $(\C^*)^{n+1}\supset\C\hat\Del$.

{\it Step 2}. Now we show that $\C Ch(\hat F)\cap H_a$,
$a\in[0,\eps)$, is an isotopy in $\C\hat\Del$.
For $\eps>0$ small enough,
$\{\hat F = 0\}\cap\{t=a\} = \{F=0\}\cap\{t=a\}$
and the intersection is transverse
in $(\C^*)^{n+1}$. This implies that the family
$\C Ch(\hat F)\cap H_a$,
$a\in(0,\eps)$, is an isotopy in $\C\hat\Del$.

Now we pick a point
$\ow=(y_1v^0_1,...,y_nv^0_n,y_{n+1})
\in H_0\cap\C Ch(\hat F)$, where
$(y_1,...,y_n)\in I(\del)$ for some $s$-dimensional
face of $\Del_1,...,
\Del_N$, $y_{n+1}=\nu(y_1,...,y_n)$,
and $(v^0_1,...,v^0_n)\in(S^1)^n$. Following the proof of
Proposition \ref{l1} and using transformations
$$(x_1,...,x_n,t)\mapsto\left(t^{\bet_1}\prod_{k=1}^nx_k^{\alp_{k1}},
...,t^{\bet_n}\prod_{k=1}^nx_k^{\alp_{kn}},t\right)$$
with $A=(\alp_{kj})\in SL(n,\Z)$, we can get
$\del$ to be parallel to the coordinate $s$-plane $\{i_1=...=0\}$,
$\nu\big|_\del=\nu_0$, and $\nu\big|_{\Del\backslash\del}>\nu_0$.
Then we introduce a parameterization of a neighborhood
$U_\ow$ of the point $\ow$ in $\C\hat\Del\cap\{\im(w_{n+1})=0\}$
by $\pr_s(U')\times U_0\times U''$ via an extension of the
map $\tmu_{\hat\Del}$ such that, in the notation
of the proof of Proposition \ref{l1}, $\pr_s(U')$ is given by
the first $s$ inequalities in (\ref{e15}), $U''$ is given by
(\ref{e14}), and $U_0$ is a neighborhood of the point
$\pr_{n-s+1}(\del)$ in $\pr_{n-s+1}(\hat\Del)$.
Here $\Int(U_0)$ is identified with the domain
$\pr_{n-s+1}(U')\subset(\R^*_+)^{n-s+1}$, given in the coordinates
$x_{s+1},...,x_n,t$ by the relations
$$M_{i_{s+1},...,i_n,i_{n+1}}=x_{s+1}^{i_{s+1}}...x_n^{i_n}
t^{i_{n+1}}<\eps,\quad(i_{s+1},...,i_{n+1})\in
\pr_{n-s+1}(\hat\Del\backslash\Del)\cap\Z^{n-s+1}\ ,$$
and $\C Ch(\hat F)\cap\Int(U_\ow)$ is parameterized in
$\pr_s(U')\times\Int(U_0)\times U''$ by equations (\ref{e17}),
where~$G$ and~$H$ are vanishing at zero smooth functions
of $x_2,...,x_s$, $v_1,...,v_n$ and $M_{i_{s+1},...,i_n,i_{n+1}}$,
$(i_{s+1},...,i_{n+1})\in
\pr_{n-s+1}(\hat\Del\backslash\Del)\cap\Z^{n-s+1}$.

Denote $\widetilde\Del=\del\times\pr_{n-s+1}(\hat\Del)$.
Clearly, $\hat\Del$ and $\widetilde\Del$ coincide in a neighborhood
of $\ow$. Moreover, The maps $\tmu_{\hat\Del}$ and
$\tmu_{\widetilde\Del}$ are connected by an isotopy on
the domain $\pr_s(U')\times\Int(U_0)\times U''$: such an isotopy
can be written explicitly by supplying the non-common
summands in the formulae for $\tmu_{\hat\Del}$ and
$\tmu_{\widetilde\Del}$ by a parameter running over $[0,1]$.
This isotopy extends up to equivariant tame isotopy on
$\pr_s(U')\times U_0\times U''$, so that replacing
$\tmu_{\hat\Del}$ by $\tmu_{\widetilde\Del}$, one replace $U_\ow$ by
another neighborhood $\widetilde U_\ow$ of $\ow$ in
$\C\widetilde\Del\cap\{\im(w_{n+1})\}$.
Note that the map $\tmu_{\widetilde\Del}$ splits on
$\pr_s(U')\times U_0\times U''$ into the product of
$$\mu_\del:\pr_s(U')\to\del,\quad\mu_{\pr_{n-s+1}(\hat\Del)}\to
\pr_{n-s+1}(\hat\Del),\quad\Id:U''\to U''.$$
Together with the result of Step 1 this allows us to introduce
in $U_0$ coordinates $t,\tet_1,...,\tet_{n-s}$ so that
$$U_0=\{0\le t<\eps,\ -\eps<\tet_j<\eps\},\quad U_0\cap
\partial(\pr_{n-s+1}(\hat\Del))=\{t=0\}\ ,$$
and $\tmu_{\widetilde\Del}(H_a)=\pr_s(U')\times\{t=a\}\times U''$.
In addition, $\mu_{\pr_{n-s+1}(\hat\Del)}$ expresses
$M_{i_{s+1},...,i_n,i_{n+1}}$,
$(i_{s+1},...,i_{n+1})\in
\pr_{n-s+1}(\hat\Del\backslash\Del)\cap\Z^{n-s+1}$, as continuous
functions of $t,\tet_1,...,\tet_{n-s}$, vanishing at zero. Hence
the closure of $\tmu_{\widetilde\Del}(\{F=0\})$ in
$\widetilde U_\ow$ is given by equations
\begin{equation}
x_1=x^0_1+G',\quad v_1=v^0_1+H',\label{e18}
\end{equation}
where~$G'$ and~$H'$ are vanishing at zero smooth functions of
$x_2,...,x_s$, $v_2,...,v_n$, $t,\tet_1,...,\tet_{n-s}$.
The variety (\ref{e18}) intersects any hypersurface
$\pr_s(U')\times\{t=a\}\times U''$, $a\in[0,\eps)$, transversally
in $\pr_s(U')\times U_0\times U''\simeq\widetilde U_\ow
\simeq U_\ow$, thereby proving that
$\C Ch(\hat F)\cap H_a$,
$a\in[0,\eps)$, is an isotopy in $\C\hat\Del$.

{\it Step 3}.
For $a\ne 0$, the projection of $\C Ch(\hat F)\cap H_a$ into
$\C\Del$ is the closure of the image of $\{F=0\}\cap\{t=a\}
\subset(\C^*)^{n+1}$ by the map
\begin{equation}
\tmu_{\Del,a}(\ox,\ov)=
\frac{\sum_{(i_1,...,i_n)\in\Del}(i_1v_1,...,i_nv_n)\cdot
\left(x_1^{i_1}...x_n^{i_n}\sum_{(i_1,...,i_n,k)\in\hat\Del}a^k
\right)}{\sum_{(i_1,...,i_n)\in\Del}
\left(x_1^{i_1}...x_n^{i_n}\sum_{(i_1,...,i_n,k)\in\hat\Del}a^k
\right)}\ .\label{e20}
\end{equation}
The map $\tmu_{\Del,a}$ is connected with $\tmu_\Del$
on $\{t=a\}\simeq(\C^*)^n$ by the isotopy
\begin{equation}
\frac{\sum_{(i_1,...,i_n)\in\Del}(i_1v_1,...,i_nv_n)\cdot
\left(x_1^{i_1}...x_n^{i_n}\left(\sum_{(i_1,...,i_n,k)\in\hat\Del}a^k
(1-\tau)+\tau\right)
\right)}{\sum_{(i_1,...,i_n)\in\Del}
\left(x_1^{i_1}...x_n^{i_n}\left(\sum_{(i_1,...,i_n,k)\in\hat\Del}a^k
(1-\tau)+\tau\right)\right)}\ ,\label{e21}
\end{equation}
providing a tame isotopy of the projection of $\C Ch(\hat F)\cap H_a$
with $\C Ch(F\big|_{t=a})\subset\C\Del$. Similarly, the
projection of $\C Ch(\hat F)\cap\C\hat\Del_k\cap\{\im(w_{n+1})=0\}$
into $\C\Del_k$, $1\le k\le N$, is the closure of the image
of $\{F^{\hat\Del_k}=0\}\cap\{t=1\}$ by the map
$\tmu_{\hat\Del_k}$ which coincides with the closure of
the image of $\{F_k=0\}$ by the map $\tmu_{\Del_k}$, {\it i.e.},
$\C Ch(F_k)$.

So, Theorem is proven in the complex case.
If all the numbers in~${\cal A}$ are real, then
the isotopy constructed is
equivariant.
\proofend

\begin{theorem}\label{t3}
{\rm (Projective complex Viro theorem, see~\cite{Vi1})}.
Let, in the notation of Theorem \ref{l5}, $\Del=T_d^n$ and the
subdivision $\cal S$ be defined by a convex piecewise-linear
function $\nu:T_d^n\to\R$. Then $P\C Ch({\cal S},{\cal A})$ is
tame isotopic in $\C P^n$ (equivariantly, if all the numbers
in~${\cal A}$ are real) to $P\C Ch(F)$, where $$F(z_1,...,z_n,t)=
\sum_{i_1+...+i_n\le d}A_{i_1...i_n}t^{\nu(i_1,...,i_n)}
z_1^{i_1}\cdot... \cdot z_n^{i_n}=0,$$ and $t=\const>0$ is
sufficiently small.
\end{theorem}

{\bf Proof.}
We multiply all the polynomials by $z_1...z_n$ and apply Theorem
\ref{l5} to the simplex $\widetilde T_d^n=T_d^n+(1,...,1)$
and correspondingly shifted ${\cal S},{\cal A}$. Then we note that
the isotopy of $\C Ch(z_1...z_nF)$ and $\bigcup_{k=1}^N
\C Ch(z_1...z_nF_k)$ in $\C\widetilde T_d^n$ is tame and
compatible with the action of $(S^1)^n$ on
$\partial\C\widetilde T_d^n$. This allows us to take
quotient by this action as was done in the proof of
Proposition \ref{l3} and obtain the required isotopy.
\proofend

\begin{proposition}\label{c3}
Let two subdivisions ${\cal S}=\{\Del_k,\ k=1,...,N\}$ and ${\cal
S}'=\{\Del_{kl},\ l=1,...,r_k,\ k=1,...,N\}$ of a polytope
$\Del\subset(\R^*_+)^n$ satisfy
$$\Del_k=\bigcup_{l=1}^{r_k}\Del_{kl},\quad k=1,...,N,$$ so that
the subdivision ${\cal S}'$ is given by piecewise-linear function
$\nu:\Del\to\R$, whose restrictions $\nu_k=\nu\big|_{\Del_k}$,
$k=1,...,N$, are convex. Then the varieties $\C Ch({\cal S},{\cal
A})$ and $\C Ch({\cal S}',{\cal A}')$ are tame isotopic, provided
${\cal A},{\cal A}':\Del\to\C$ define non-degenerate polynomials
$F_k$, $F_{kl}$, $l=1,...,r_k$, $k=1,...,N$. Similarly, given two
subdivisions $\cal S$, ${\cal S}'$ of $T_d^n$ satisfying the
previous assumptions, the C-hypersurfaces of degree~$d$
constructed out of ${\cal S},{\cal A}$ and ${\cal S}',{\cal A}'$
are tame isotopic in $\C P^n$.
\end{proposition}

{\bf Proof.} Consider the case $\Del\subset(\R^*_+)^n$.
Fix $a>0$. Let the maps $\tmu_{\Del_k,a}:(\C^*)^n\to\C\Del_k$,
$k=1,...,N$, be defined by (\ref{e20}), where summations run
over $\Del_k$ and $\hat\Del_k$, constructed as in the proof of
Theorem \ref{l5} by means of $\nu_k$.
Put
$$F_{k,a}=\sum_{(i_1,...,i_n)\in\Del_k}A'_{i_1...i_n}z_1^{i_1}...
z_n^{i_n}a^{\nu(i_1,...,i_n)},\quad k=1,...,N.$$
As shown in the proof of Theorem \ref{l5}, the closure
$\C Ch_a(F_{k,a})$ of $\tmu_{\Del_k,a}(\{F_{k,a}=0\})$ in
$\C\Del_k$ is tame isotopic to $\bigcup_l\C Ch(F_{kl})$.
Moreover, for any face $\del=\Del_k\cap\Del_j$,
the isotopies in~$\C\Del_k$ and~$\C\Del_l$ coincide on~$\C\del$.
Then, for each $k=1,...,N$, we connect $\tmu_{\Del_k,a}$
with $\tmu_{\Del_k}$ by an isotopy (\ref{e21}),
thereby obtaining the required isotopy of
$\C Ch({\cal S},{\cal A})$ and
$\C Ch({\cal S}',{\cal A}')$.

The same argument proves the statement for $\Del=T_d^n$.
\proofend

A subdivision $\Del=\Del_1\cup...\cup\Del_N$ is called
{\it maximal} if it cannot be refined. In this case all the
integral points in $\Del$ are vertices of $\Del_1,...,\Del_N$,
and these polytopes are simplices.

\begin{corollary}\label{c5}
Given a polytope $\Del\subset(\R^*_+)^n$, any C-hypersurface in
$\C\Del$ is tame isotopic to a C-hypersurface constructed out of
a maximal subdivision of $\Del$. The same is true for C-hypersurfaces
in $\C P^n$.
\end{corollary}

{\bf Proof.} Let $\Del=\Del_1\cup...\cup\Del_N$, and $f:\Del\to\R$
be a smooth convex function. Define a piecewise-linear function
$\nu:\Del\to\R$ as follows: put $\nu(\oi)=f(\oi)$,
$\oi\in\Del\cap\Z^n$, then define the graph of $\nu\big|_{\Del_k}$
to be the ``lower" part of the convex hull of $\{(\oi,\nu(\oi)\ :\
\oi\in\Del_k\}$. The function $\nu$ defines a maximal subdivision
of $\Del$ inscribed into the initial subdivision and satisfying
the conditions of Proposition \ref{c3}, which completes the proof.
\proofend

\begin{corollary}\label{c4}
Any C-curve in $\C P^2$ is tame isotopic
(not equivariantly, in general) to an
algebraic curve of the same degree.
\end{corollary}

{\bf Proof.}
By Corollary \ref{c5} we can assume that a C-curve in
$\C P^2$ is constructed out of a maximal triangulation of
$T_d^2$. We will transform any given triangulation into
a convex triangulation,
so that in each transformation step the conditions
of Proposition \ref{c3} hold true.

Let $\cal S$ be a triangulation of $T_d^2$, $O({\cal S})$
denote the star of the origin $O$ with respect to this triangulation,
and outside $O({\cal S})$ the triangulation $\cal S$ is
maximal.
We construct a triangulation ${\cal S}'$ of $T_d^2$ such that
$O({\cal S}')\supset O({\cal S})$ and $O({\cal S}')\ne O({\cal S})$.
Then, in finitely many steps we come to $O({\cal S})=
T_d^2$, which corresponds to a convex triangulation,
and the required statement will follow from Theorem \ref{t3} and
Proposition \ref{c3}.

Let $O,P_1,...,P_r$ be all the vertices of $O({\cal S})$ numbered
successively clockwise along
$\partial O({\cal S})$. Any segment $[P_i,P_{i+1}]$
either lies on $\partial T_d^2$, or is an edge of a
unique triangle
$T_i\in{\cal S}$, $T_i\not\subset O({\cal S})$.

\begin{figure}
\begin{center}
\epsfxsize 145mm
\epsfbox{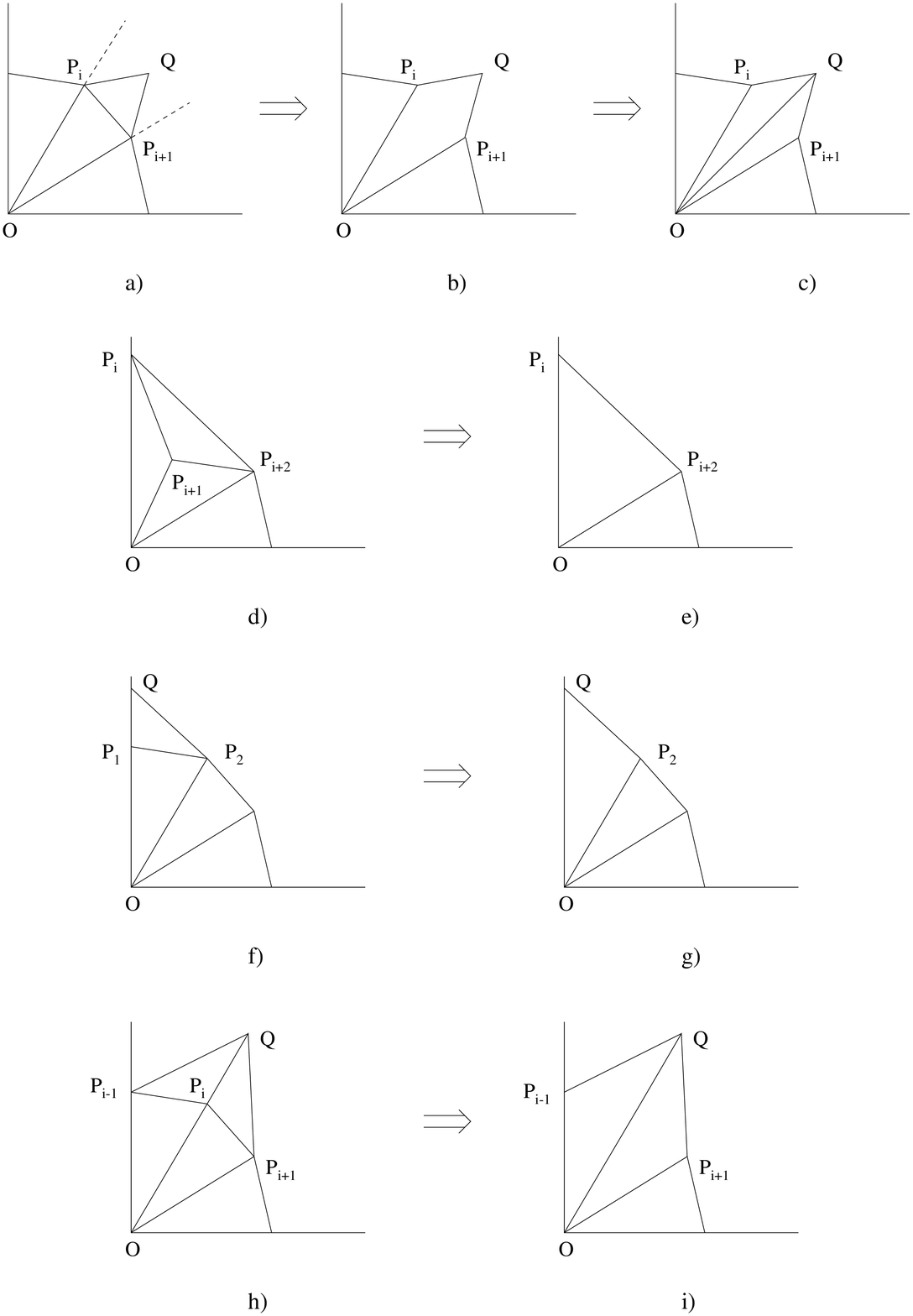}
\end{center}
\caption{Transformation of a triangulation}
\label{f1}
\end{figure}

(i) Assume that, for some $i=1,...,r$,
the vertex $Q\ne P_i,P_{i+1}$ of $T_i$ lies
between the straight lines$(OP_i)$ and
$(OP_{i+1})$, and $Q\not\in(OP_i)\cup(OP_{i+1})$
(see Figure \ref{f1}a). Then we change the subdivision of $T_d^2$
as shown in Figure \ref{f1}b,c. These changes satisfy the conditions
of Proposition \ref{c3} and lead to a triangulation with a strongly
greater star of $O$.

(ii) Assume that, for some $i=1,...,r-1$,
$T_i=T_{i+1}$, {\it i.e.}, the vertices of the latest triangle are
$P_i,P_{i+1},P_{i+2}$ (see Figure \ref{f1}d).
Then we perform the transformation shown in Figure
\ref{f1}e, once again
increasing the star of $O$.

(iii) Assume that there are no triangles $T_i$ as in (i), (ii).
Then any triangle $T_i$ is either ``left", {\it i.e.}, the vertex
$Q$ lies on $(OP_i)$ or above $(OP_i)$, or ``right",
{\it i.e.}, $Q$ lies on $(OP_{i+1})$ or below $(OP_{i+1})$.
If there exist ``left" triangles, consider the ``left" triangle
$T_i$ with the minimal $i$. If $i=1$, we have the situation shown
in Figure \ref{f1}f. Then we change triangulation as shown
in Figure \ref{f1}g, increasing the star of $O$. If $i>1$, then
the triangle $T_{i-1}$ must be ``right", which means that we have
a situation shown in Figure~\ref{f1}h. Then
we change triangulation as shown
in Figure \ref{f1}i, increasing the star of $O$.
\proofend

\subsection{Basic properties of C-hypersurfaces}\label{first-prop}

The real part of a real C-hypersurface~$M$ in~$\C P^n$
(see Definition~\ref{d1})
is denoted
by~$\R M$.

\begin{proposition}\label{l6}
A complex C-hypersurface~$M$ of degree~$d$ in~$\C P^n$
is an orientable manifold,
homologous to an algebraic
hypersurface of degree~$d$ in~$\C P^n$.
If~$M$ is real its real part~$\R M$
is a closed
manifold, mod 2 homologous to the real point set of a
real algebraic hypersurface of degree~$d$
in~$\R P^n$.
\end{proposition}

{\bf Proof.}
The Jacobian of the moment map
$\mu_{\Del}:(\R^*_+)^n\to\Int(\Del)$,
$\mu_\Del(\ox)=(\mu_\Del^{(1)}(\ox),...,\mu_\Del^{(n)}(\ox))$,
is positive. Hence the Jacobian of the extended moment map
$\tmu_\Del:(\C^*)^n\to\C I(\Del)$ is positive in the coordinates
$$x_1=|z_1|,\ v_1=\frac{z_1}{|z_1|},\ ...\ ,\ x_n=|z_n|,\ v_n
\frac{z_n}{|z_n|}\ ,$$
because this is a diffeomorphism, and at a point with
$v_1=...=v_n=1$ one can easily compute
$$\det\left(\frac{D(\tmu_{\Del})}{D(x_1,v_1,...,x_n,v_n)}\right)=
\det\left(\frac{D(\mu_{\Del})}{D(x_1,...,x_n)}\right)\cdot\prod_{j=1}^n
\mu_{\Del}^{(j)}(x_1,...,x_n)>0\ .$$

This means, in particular, that $\tmu_{\Del}$ canonically defines
an orientation of images of complex subvarieties of
$(\C^*)^n$. Therefore, the open subset
$$\bigcup_{k=1}^N(\C Ch(F_k)\cap\C I(\Del_k))$$
of $\C Ch({\cal S},{\cal A})$ is canonically orientable.
So, to complete the proof of orientability, one should verify that
these orientations are compatible when gluing the charts $\C Ch(F_k)$,
$\C Ch(F_j)$ with $\Del_k\cap\Del_j=\del$ being a common facet
(the gluing along faces of lower dimensions does not affect the
orientation).
In other words,
the orientations of $\C Ch(F_k^{\del})=\C Ch(F_j^{\del})$
induced by $\C Ch(F_k)$ and $\C Ch(F_j)$ are opposite.

Without loss of generality suppose that $\del$ is contained in a
hyperplane $i_1=a>0$, $\Del_k\subset\{i_1\le a\}$, and
$\Del_j\subset\{i_1\ge a\}$. Then by construction, $\C\Del_k$ induces
on $\C\del$ an orientation defined by the form
$$-dv_1\wedge dx_2\wedge dv_2\wedge...\wedge dx_n\wedge dv_n,$$
and $\C\Del_j$ induces on $\C\del$ the opposite orientation defined by
$$dv_1\wedge dx_2\wedge dv_2\wedge...\wedge dx_n\wedge dv_n.$$

On the other hand, a coorienting 2-vector bundle on $\C Ch(F_k)\cap
\C I(\Del_k)$ can continuously be extended to a coorienting
2-bundle on $\C Ch(F_k^{\del})\subset\C I(\del)$, and the same
for $\Del_j$.
Indeed, let the complex straight line $\Lam$
$$z_2=...=z_n=xv,\quad z_1=x_0=\const,\quad x\in(0,\infty),\ |v|=1,$$
with small $z_0>0$ meet the hypersurfaces
$\{F_k=0\}$,
$\{F_k^{\del}=0\}$ transversally. Then the family of surfaces
$\Sig_{\lam}$, $\lam\in[0,1]$,
$$\left\{\frac{\sum_{(i_1,...,i_n)\in\Del_k}x^{i_2+...+i_n}
x_0^{i_1}\lam^{a-i_1}\cdot(i_1,i_2v,...,i_nv)}
{\sum_{(i_1,...,i_n)\in\Del_k}x^{i_2+...+i_n}
x_0^{i_1}\lam^{a-i_1}}\ :\ x\in(0,\infty),\ |v|=1\right\}$$
is a diffeotopy,
connecting $\Sig_0=\tmu_{\del}(\Lam)$, which coorients
$\C Ch(F_k^{\del})$ in $\C I(\del)$, and $\Sig_1=
\tmu_{\Del_k}(\Lam)$, which coorients $\C Ch(F_k)$ in $\C I(\Del_k)$.

Comparing this with the previous remark on orientations of $\C\del$,
one completes the proof of the orientability of a complex
C-hypersurface.

At last,
$[\C Ch({\cal S},{\cal A})]\in
H_{2n-2}(\C P^n)$ and
$[\R Ch({\cal S},{\cal A})]\in H_{n-1}(\R P^n,\Z/2)$
can easily be computed by induction considering the intersection
of the complex and real charts with the coordinate
planes.
\proofend

\begin{proposition}\label{l7}
Any (real) C-hypersurface $M$ is (equivariantly)
tame isotopic to a close smooth
manifold $M_{sm}$ of codimension 2 in $\C P^n$.
\end{proposition}

{\bf Proof.} Let $M=\C Ch({\cal S},{\cal A})/S^1$, where
${\cal S}=\{\Del_1,...,\Del_N\}$, ${\cal A}=
\{A_\oi\ :\ \oi\in T_d^n\cap\Z^n\}$.
Put $F_m = \sum_{\oi\in\Del_m}A_\oi\oz^\oi$, $m=1,...,N$.
We construct two nonvanishing
$\R$-linearly independent sections~$\os$ and~$\os'$ of the
bundle $T\C P^n\big|_{U(M)}$, where $U(M)$ is a neighborhood of
$M$ in $\C P^n$, such that $\os$ is
equivariant, $\os'$ is anti-equivariant, and
the 2-subbundle
$\Span_\R\{\os,\os'\}\subset T\C P^n\big|_M$
is transverse to $T\C Ch(F_m)$ in $T\C P^n\big|_{\C Ch(F_m)}$
for any $m=1,...,N$. This will imply
the existence of a smooth $2(n-1)$-manifold $M_{sm}$ isotopic
(equivariantly, if $\conj(M)=M$) to $M$, close to $M$ and transverse
to the 2-bundle $\Span_\R\{\os,\os'\}$. Namely, one smoothes
$M$, pushing it along the trajectories of the
vector field $\os$ in a neighborhood of
$\bigcup_{m=1}^N\partial\C\Del_m$.

(i) First, we shift $T_d^n$ into $\widetilde T_d^n=(1,...,1)+
T_d^n$ (shifting respectively ${\cal S},{\cal A}$ as well), and
construct sections~$\os$ and~$\os'$ of~$T\C^n$ defined in
a neighborhood of $\C Ch({\cal S},
{\cal A})$, with the above properties and an additional one:
invariance with respect to the $S^1$-action on
$\partial\C\widetilde T_d^n$. The latter allows us to
obtain the required sections of $T\C P^n$ along the procedure
described in the proof of Proposition
\ref{l3}.

(ii) Fix $m=1,...,N$ and consider the hypersurface
$\widetilde F_m(\oz)=F_m(e^{z_1},...,e^{z_n})=0$
in $\C^n$. Define vector fields~$\os_1$ and~$\os'_1$ on~$\C^n$ by
$$\os_1(\oz)=\frac{\conj(\grad\widetilde F_m)}
{|\conj(\grad\widetilde F_m)|},\quad
\os'_1(\oz)=\frac{\conj(\grad\widetilde F_m)\sqrt{-1}}
{|\conj(\grad\widetilde F_m)\sqrt{-1}|}\ .$$
They do not vanish along
$\{\widetilde F_m=0\}$ and span a 2-bundle
orthogonal to
$T\{\widetilde F_m=0\}$ in $T\C^n\big|_{\{\widetilde F_m=0\}}$.
These vector fields are $2\pi\sqrt{-1}$-periodic in
each coordinate $z_1,...,z_n$, and
their normalized images~$\os$ and~$\os'$ by the differential
$D(\tmu_{\Del_m}\circ\exp)$ of the map
$\tmu_{\Del_m}\circ\exp:\C^n\to\C\Del_m$ give a 2-bundle
$\Span_\R\{\os,\os'\}$ on $\C I(\Del_m)\cap\C Ch(F_m)$ transverse to
the tangent bundle $T(\C I(\Del_m)\cap\C Ch(F_m))$. We will
show that these vector fields continuously extend to
$\bigcup_{m=1}^N(\partial\C\Del_m)\cap U(M)$.

(iii) Let~$\del$
be a face of $\Del_m$, $\ow$ be a point in $\C Ch(F_m)\cap\C I(\del)$,
and $U_\ow\subset\C\Del_m$ be a neighborhood of $\ow$ as
introduced in the step (iv) of the proof of Proposition \ref{l1}.
We claim that, for any point $\oz\in\C^n$ such that
$\tmu_\Del(e^\oz)\in U_\ow\cap\C I(\Del_m)$,
one has
\begin{equation}
\bigg|\bigg|\frac{D(\tmu_{\Del_m}\circ\exp)_\oz}
{||D(\tmu_{\Del_m}\circ\exp)_\oz||}-
\frac{D(\tmu_\del\circ\exp)_{\overline\alp}}
{||D(\tmu_\del\circ\exp)_{\overline\alp}||}
\bigg|\bigg|<c\eps\ ,\label{e35}
\end{equation}
and
\begin{equation}
\bigg|\frac{\grad(F_m\circ\exp)}{|\grad(F_m\circ\exp)|}(\oz)
-\frac{\grad(F_m^\del\circ\exp)}{|\grad(F_m^\del\circ\exp)|}
(\overline\alp)
\bigg|<c\eps\ ,\label{e37}
\end{equation}
where $\overline\alp\in\C^n$ is some point,
satisfying $\tmu_\Del(\exp(\overline\alp))=\ow$, and
$c>0$ depends only on the coefficients of $F$ and
on the point $\ow$.

(iv) To show (\ref{e35}), note that $\tmu_{\Del_m}\circ\exp$ splits
into $\mu_{\Del_m}\circ\exp:\R^n\to I(\Del_m)$ and
$\exp:\R^n\sqrt{-1}\to(S^1)^n$. Inequality (\ref{e35})
for the differential of $\exp:\R^n\sqrt{-1}\to(S^1)^n$
instead of $\tmu_{\Del_m}\circ\exp$ immediately follows from
(\ref{e14}). Then, acting by a transformation from
$SL(n,\Z)$ and a shift, we move $\Del_m,\del$ into
$\widetilde\Del_m,\widetilde\del$ such that $\widetilde\del$
lies in a coordinate $s$-plane. According to
(\ref{e15}),
$$|e^{(\oi,\real\oz)}-
e^{(\oi,\real\overline\alp)}|<c_1\eps,\
\oi\in\widetilde\del\cap\Z^n,\quad
e^{(\oi,\real\oz)}<c_1\eps,\ \oi\in(\widetilde\Del_m
\backslash\widetilde\del)\cap\Z^n,$$
where $c_1>0$ depends only on $\Del_m$, thereby this
implies (\ref{e35}),
because $D(\mu_{\widetilde\Del_m}\circ\exp)$ is represented by
the matrix
$$\left(\frac{\sum_{\oi\in\widetilde\Del_m}i_pi_qe^{(\oi,\real\oz)}
\cdot\sum_{\oi\in\widetilde\Del_m}e^{(\oi,\real\oz)}-
\sum_{\oi\in\widetilde\Del_m}i_pe^{(\oi,\real\oz)}\cdot
\sum_{\oi\in\widetilde\Del_m}i_qe^{(\oi,\real\oz)}}
{(\sum_{\oi\in\widetilde\Del_m}e^{(\oi,\real\oz)})^2}\right)_{p,q
=1,...,n}\ ,$$
and $D(\mu_{\widetilde\del}\circ\exp)$ does not vanish.

(v) Similarly,
$F_m(e^\oz)=e^{(\oi_0,\oz)}\sum_{\oi\in\Del}
A_\oi e^{(\oi,\oz)}$,
where by~(\ref{e14}) and~(\ref{e15})
$$|e^{(\oi,\oz)}-e^{(\oi,\overline\alp)}|<c_2\eps,\ \oi\in\del
\cap\Z^n,
\quad |e^{(\oi,\oz)}|<c_2\eps,\ \oi\in(\Del\backslash\del)
\cap\Z^n,$$
where $f_2>0$ depends on $\Del_m$.
So, we obtain
$$\grad(F_m\circ\exp)(\oz)=e^{(\oi_0,\oz)}\left(\sum_{\oi\in\del}
A_\oi\oi e^{(\oi,\oz)}+\sum_{\oi\in\Del\backslash\del}
A_\oi\oi e^{(\oi,\oz)}\right)\ ,$$
which immediately implies (\ref{e37}), because
$\sum_{\oi\in\del}A_\oi\oi e^{(\oi,\overline\alp)}
=\grad(F_m^\del\circ\exp)(\overline\alp)\ne 0$.

(vi) Relations~(\ref{e35})
and~(\ref{e37}) provide continuous extension
of the vector fields~$\os$ and~$\os'$
on $\bigcup_{m=1}^N(\partial\C\Del_m)\cap U(M)$ so that
they remain $\R$-linearly independent and belong to
$T\C I(\del)$ along $\C I(\del)$ for any proper face
$\del$ of $\Del_m$, $m=1,...,N$. Moreover, the restrictions
of~$\os$ and~$\os'$ on $\C I(\del)$ depend only on~$\del$ and
$F_m^\del$, common for all $\Del_m\supset\del$, hence they are
compatible with $S^1$-action on $\partial\C\widetilde T_d^n$,
and we are done.
\proofend

\begin{proposition}\label{p9}
Given a (real) C-hypersurface $M$,
the tangent bundle
to its smoothing $M_{sm}$
is (equivariantly) isotopic to a
(equivariant) bundle of complex hyperplanes.
In particular, $M_{sm}$ possesses an (equivariant) almost complex structure.
If~$M$ is real, the above isotopy has fixed
intersection with $T\R P^n\big|_{M_{sm}}$ (equal to $T\R M_{sm}$).
Furthermore, the complex structure in $\C P^n$ can be
(equivariantly) deformed into an almost complex structure,
compatible with the metric and
for which $TM_{sm}$ is invariant.
\end{proposition}

{\bf Proof.} We construct the required isotopy in few steps.

{\it Step 1}.
In the proof of Proposition \ref{l7} we have constructed
sections $\os$, $\os'$ of the bundle $T\C P^n\big|_{M_{sm}}$,
which are linearly
independent at any point $\ow\in M_{sm}$ and such that the 2-bundle
$\Span_\R\{\os,\os'\}$ is
transversal to $TM_{sm}$. Suppose that $\C P^n$ is equipped with
a Hermitian metric compatible with the complex structure and
the complex conjugation. Then there exists a (linear) isotopy of the
bundle $TM_{sm}$ into the $(2n-2)$-bundle
$\left(\Span_\R\{s,s'\}\right)^\perp$.

{\it Step 2}. Assume that $M=P\C Ch({\cal S},{\cal A})$, where $\cal S$
is a subdivision $T_d^n=\Del_1\cup...\cup\Del_N$.
Via the map $\nu_d^n:\C T_d^n\to\C P^n$ we
pull back $\os$ and $\os'$ to sections of the bundle
$T\C T_d^n$ restricted to $(\nu_d^n)^{-1}(M_{sm})$
(which will be denoted by $M_{sm}$ for abuse of notations),
as well as pull back
the complex structure $J$ and the Hermitian metric.
Note also that $T\C T_d^n\simeq\C T_d^n\times\C^n$ possesses
the standard complex structure being just the multiplication
by $\sqrt{-1}$.

All this data on $\C T_d^n$ is
invariant with respect to the $S^1$-action on
$\partial\C T_d^n$, as well as the procedures used further; hence
the isotopies we construct in $T\C T_d^n$ can be pushed to
$T\C P^n$.

The following lemma implies that there exists an isotopy of the bundle
$\left(\Span_\R\{\os,\os'\}\right)^\perp$ into the bundle
$\left(\Span_\R\{\os,\os\sqrt{-1}\}\right)^\perp$.

\begin{lemma}\label{l12}
In the above notation, for any $\ow\in M_{sm}$, no vector
$(1-\lam)\os'(\ow)+\lam\os\sqrt{-1}$, $0\le\lam\le 1$, is
$\R$-proportional
to $\os(\ow)$.
\end{lemma}

{\bf Proof.} Let $\ow\in\C\Del_m$, $1\le m\le N$.
We have
$$\os(\ow)=D\widetilde\tmu_{\Del_m}(\ob),\quad
\os'(\ow)=D\widetilde\tmu_{\Del_m}(\ob\sqrt{-1}),$$
where $\widetilde\tmu_{\Del_m}:\C^n\to\C^n$ is defined by
\begin{equation}
\C^n=\R^n\oplus\R^n\sqrt{-1}\stackrel{\log
\widetilde\mu_{\Del_m}\oplus\Id}
{\longrightarrow}\R^n\oplus\R^n\sqrt{-1}=\C^n\ ,\label{e36}
\end{equation}
$$\log
\widetilde\mu_{\Del_m}=(\log\mu_{1,\Del_m},...,\log\mu_{n,\Del_m}),
\quad\mu_{\Del_m}=(\mu_{1,\Del_m},...,\mu_{n,\Del_m})\ ,$$
and $\ob$ is a nonzero vector on $\C^n$.

Assume that
\begin{equation}
\kappa\os(\ow)=
(1-\lam)\os'(\ow)+\lam\os(\ow)\sqrt{-1}\label{e100}
\end{equation}
for some $\kappa\in\R$,
$\lam\in(0,1)$ and $\ow\in\C I(\Del_m)$.
Splitting $\ob$ into $\ob_r+\ob_i\sqrt{-1}$
with $\ob_r,\ob_i\in\R^n$, we obtain in view of (\ref{e36})
\begin{eqnarray}
&\os(\ow)=D(\log\mu_{\Del_m})(\ob_r)+\ob_i\sqrt{-1},\nonumber\cr
&\os(\ow)\sqrt{-1}
=-\ob_i+D(\log\mu_{\Del_m})(\ob_r)\sqrt{-1},\nonumber\cr
&\os'(\ow)=-D(\log\mu_{\Del_m})(\ob_i)+\ob_r\sqrt{-1}\ .\nonumber
\end{eqnarray}

Suppose that $\kappa=0$. Setting the expressions for $\os(\ow)\sqrt{-1}$
and $\os'(\ow)$ into (\ref{e100}) we get that the operator
$D(\log\mu_{\Del_m})$ has a negative eigenvalue, what is impossible.
Indeed, the matrix of $D(\log\mu_{\Del_m})$ is a
product of
$$A=\left(\frac{\sum_{\oi\in\Del}i_pi_qe^{(\oi,\ox)}\cdot
\sum_{\oi\in\Del}e^{(\oi,\ox)}-(\sum_{\oi\in\Del}i_pe^{(\oi,\ox)})
(\sum_{\oi\in\Del}i_qe^{(\oi,\ox)})}{(\sum_{\oi\in\Del}e^{(\oi,
\ox)})^2}\right)_{p,q=1,...,n}\ ,$$
$$B=\diag\left(\frac{\sum_{\oi\in\Del}i_pe^{(\oi,\ox)}}
{\sum_{\oi\in\Del}e^{(\oi,\ox)}}\right)_{p=1,...,n}\ .$$
Here both $A$ and $B$ are positive definite, for instance,
$$(A\ob,\ob)=\frac{\sum_{\oi,\oj\in\Del_m}e^{(\oi,\ox)}
e^{(\oj,\ox)}(\sum_{p=1}^n(i_p-j_p)b_p)^2}
{(\sum_{\oi\in\Del_m}e^{(\oi,\ox)})^2}>0\ ,$$
hence $AB$ rotates any vector by an angle $<\pi$, so cannot have
negative eigenvalues.

Suppose that $\kappa\ne 0$. Plugging the above expressions for
$\os(\ow)$, $\os(\ow)\sqrt{-1}$ and $\os'(\ow)$ into (\ref{e100}),
we obtain that the operator
$$(\lam-\lam^2)\Id+(1-2\lam+2\lam^2+\kappa^2)D(\log\mu_{\Del_m})+
(\lam-\lam^2)D(\log\mu_{\Del_m})^2$$
vanishes at $b_r\ne 0$. Both the roots of the polynomial
$$\varphi(X)=\lam-\lam^2+(1-2\lam+2\lam^2+\kappa^2)X+(\lam-\lam^2)X^2$$
are negative; hence $D(\log\mu_{\Del_m})$ should have a negative
eigenvalue in contrary to the previous argument, and we are done.

The same argument proves
the required statement when
$\ow \in\partial\C\Del_m$
in view of~(\ref{e35}) and~(\ref{e37}).
\proofend

{\it Step 3}.
Now we claim that there exists an isotopy of the bundle
$\left(\Span_\R\{\os,\os\sqrt{-1}\}\right)^\perp$
into the bundle
$\left(\Span_\R\{\os,J\os\}\right)^\perp$ which is a bundle of
complex hyperplanes. We use the fact that, for any $\ow\in
\C T_d^n$ and $\lam\in(0,1)$, no vector $\lam\os(\ow)\sqrt{-1}+(1-\lam)
J\os(\ow)$ is $\R$-proportional to $\os(\ow)$, which follows from
Lemma \ref{l12} applied to the case $\Del_m=T_d^n$.

{\it Step 4}. If $M$ is real then the initial 2-bundle
$\Span_\R\{\os,\os'\}$ and all the constructions are $\conj$-invariant;
hence the isotopies and the resulting bundle of complex hyperplanes are
equivariant.
To satisfy the condition that
the isotopy has fixed
intersection with $T\R P^n\big|_{M_{sm}}$ equal to $T\R M_{sm}$,
in the very beginning we choose a hermitian metric on
$\C P^n$, compatible with the complex structure and
such that the vector field $\os\big|_{\R M_{sm}}
\subset T\R P^n$ is orthogonal to $T\R M_{sm}$.

{\it Step 5}. Let us extend the above constructed isotopy to
a neighborhood $U$ of $M_{sm}$. Take a smooth function $\rho:
\C P^n\to[0,1]$ (equivariant in the real case), which is $0$ in
$\C P^n\backslash U$ and $1$ on $M_{sm}$. We have three $(2n-2)$-bundles on
$U$:
$$L^{(0)}=\left(\Span\{\os,\os'\}\right)^\perp,\quad
L^{(1/2)}=\left(\Span\{\os,\os\sqrt{-1}\}\right)^\perp,\quad
L^{(1)}=\left(\Span\{\os,J\os\}\right)^\perp\ .$$
Let $L^{(t)}$ be an isotopy of $L^{(0)}$ to $L^{(1)}$ via
$L^{(1/2)}$. By constructions in Steps 2, 3, there exist
families of isometries
$$Q_\lam:T\C P^n\big|_U\to
T\C P^n\big|_U,\ 0\le\lam\le\frac{1}{2},\quad Q_\lam(L^{(\lam)})=L^{(1/2)},
\quad Q_{1/2}=\Id,$$
$$R_\lam:T\C P^n\big|_U\to
T\C P^n\big|_U,\ \frac{1}{2}\le\lam\le 1,\quad R_\lam(L^{(\lam)})=L^{(1)},
\quad R_1=\Id\ .$$
Then we define a deformation $J_t$, $t\in[0,1]$,
of the complex structure $J=J_0$
by
$$J_t=R_{1-t\rho(\ow)}J(R_{1-t\rho(\ow)})^{-1},
\quad t\rho(\ow)\le\frac{1}{2},$$
$$J_t=Q_{1-t\rho(\ow)}R_{1/2}J(Q_{1-t\rho(\ow)}R_{1/2})^{-1},
\quad t\rho(\ow)>\frac{1}{2}\ . \quad \text{\proofend}$$

\subsection{Algebraic covering}\label{alg-cover}

Introduce the map $\Pi_m:\C P^n\to\C P^n$,
$[z_0: \ldots : z_n] \mapsto
[z_0^m: \ldots :z_n^m]$.

\begin{proposition}\label{l8}
For any $C$-hypersurface~$M$ of degree $d$ in $\C P^n$
there exists a number $m_0>0$ such that for any integer $m > m_0$
the preimage
$\Pi_m^{-1}(M)$ of~$M$ is a PL-manifold tame isotopic
to a smooth algebraic
hypersurface of degree $md$ in $\C P^n$.
\end{proposition}

{\bf Proof.}
Let $M$ be
defined by polynomials $F_i(z_1,...,z_n)$, $i=1,...,N$, with
Newton polytopes $\Del_1,...,\Del_N$, where
$\Del_1\cup...\cup\Del_N= T_d^n$. The polynomials
$$\widetilde F_i(z_1,...,z_n)=
F_i(z_1^m,...,z_n^m),\quad i=1,...,N,$$
have the Newton polytopes $m\Del_i$,
$i=1,...,N$, and
define a $C$-hypersurface $\widetilde M$ of degree
$md$.

\begin{lemma}\label{p5}
The manifold $\widetilde M$ is tame isotopic in $\C P^n$ to
$\Pi_m^{-1}(M)$.
\end{lemma}

{\bf Proof.} Let $\Del\subset\R^n$ be a polytope of dimension~$n$.
Let us show that
the diffeomorphism $\mu_\Del
\Pi_m\mu_{m\Del}^{-1}:I(m\Del) \to I(\Del)$,
where $\Pi_m:(\R_+)^n\to(\R_+)^n$ acts as $(x_1,...,x_n)\mapsto
(x_1^m,...,x_n^m)$,
extends to a homeomorphism $\varphi:m\Del\to\Del$ such that,
for any proper face $\sig\subset\Del$,
the map $\varphi|_{m\sig}:m\sig\to\sig$ is a homeomorphism depending
only on $\sig$ (and not on~$\Del$).
Indeed, since $\Pi_m$ commutes with $SL(n,\Z)$ acting on
$(\R_+)^n$,
we can replace $\Del$ by its image under some $g \in SL(n,\Z)$
such that the given face $\sig$ will lie on the coordinate hyperplane
$\{i_1=0,\ i_2\cdot...\cdot i_n\ne0\}$.
Then $\mu_{t(\Del)}$ and $\mu_{mt(\Del)}$ extend
by the formulas similar to (\ref{e1}) to homeomorphisms
of $(\R_+)^n \cup\{x_1=0,\ x_2\cdot...\cdot x_n\ne 0\}$
and $I(t(\Del))\cup I(t(\sig))$ and
$I(mt(\Del))\cup I(mt(\sig))$, respectively.
This gives us a homeomorphism
$I(m\Del)\cup I(m\sig)\to I(\Del)\cup I(\sig)$,
which does not depend on the choice of $t$,
and we are done.

Consider now the following commutative diagram:
\begin{equation}
\begin{CD}
\coprod_{1\le k\le N}(\C^*)^n
@>{\coprod_{1\le k\le N}\C\mu_{m\Del_k}}>>\C T_{md}^n@>\C\mu>>\C P^n\\
@VV\coprod_{1\le k\le N}\Pi_mV @VV\varphi_mV @VV\widetilde\Pi_mV\\
\coprod_{1\le k\le N}(\C^*)^n
@>{\coprod_{1\le k\le N}\C\mu_{\Del_k}}>>\C T_d^n@>\C\mu>>\C P^n\\
\end{CD}
\nonumber\end{equation}
where $\coprod_{1\le k\le N}(\C^*)^n$ means
the disjoint union of $N$ copies of
$(\C^*)^n$.
The map $\varphi_m$ is defined on any $\C\Del_k$ as follows:
$$(x_1,...,x_n)\in\Del_k\mapsto\mu_{\Del_k}\Pi_m\mu_{m\Del_k}^{-1}
(x_1,...,x_n)\in\Del_k,$$
$$(v_1,...,v_n)\in
(S^1)^n\mapsto(v_1^m,...,v_n^m)\in(S^1)^n\ ,$$
with $\mu_{\Del_k}\Pi_m\mu_{m\Del_k}^{-1}$ extended on the whole
$\Del_k$. This definition is correct since the extensions
coming from $\Del_k$ and $\Del_j$ with a common face are the same
on the common face,
as shown above. The map $\widetilde\Pi_m$ is defined by this diagram.
Let us show that $\widetilde\Pi_m$ is tame isotopic to
$\Pi_m$. First,
$\Pi_m(x_0v_0,...,x_nv_n) = (x_0^mv_0^m,...,x_n^mv_n^m)$
is tame isotopic to
$\pi'_m(x_0v_0,...,x_nv_n)=(x_0v_0^m,
...,x_nv_n^m)$.
On the other hand, the homeomorphism
$$\coprod_{1\le k\le N}\mu_{\Del_k}\Pi_m\mu_{m\Del_k}^{-1}:
T_{md}^n\to T_d^n$$
in the definition of $\varphi_m$ is tame isotopic to
the homothety, turning $\varphi_m$ into
$$\varphi'_m(x_1v_1,...,x_nv_n)=
\frac{1}{m}(x_1v_1^m,...,x_nv_n^m)\ ,$$
and completing the proof of Lemma, since
$\pi'_m=\C\mu\circ\varphi'_m\circ(\C\mu)^{-1}$, and
$\varphi_m^{-1}(\C Ch(F_k))=\C Ch(F_k(z_1^m,...,z_n^m))$.
\proofend

\begin{lemma}\label{p6}
There exists $m_0$ such that, for any $m\ge m_0$, the subdivision
$T_{md}^n=m\Del_1\cup...\cup m\Del_N$ admits a convex refinement.
\end{lemma}

{\bf Proof.}
Let $\Gam$ be the graph of a smooth convex function
of $n$ variables, say
$f(i_1,...,i_n)=i_1^2+...+i_n^2$,
and let $\pr:\Gam\to\R^n$ be the projection.
Denote
by
$\sk^{n - 1}(\Del)$ the ($n - 1$)-skeleton
of
the subdivision $\Del_1\cup...\cup \Del_N$.
Clearly,
$\pr^{-1}(\sk^{n - 1}(\Del))$
lies on the boundary of its convex hull.
The same is true for
\begin{equation}
\pr^{-1}\left(\sk^{n - 1}(\Del)\cap\frac{1}{m}\Z^n\right). \label{e9}
\end{equation}
If $m$ is big enough, one can define a required refinement by the
piecewise-linear convex function whose graph is the lower part of
the boundary of the convex hull of the set (\ref{e9}). \proofend

Now to finish
the proof of Proposition~\ref{l8} it remains to apply
Theorem~\ref{t3} and Proposition~\ref{c3}.
\proofend

Denote by~$\chi(X)$ the Euler characteristic of~$X$, by $\sign(X)$
the signature of~$X$ if~$X$ is a manifold whose dimension
is divisible by~$4$,
and by $\chi^n_d$ (resp., $\sign^n_d$, if~$n$ is odd)
the Euler characteristic (resp., the signature) of
a nonsingular algebraic
hypersurface of degree~$d$ in~$\C P^n$.

\begin{corollary}\label{c2}
Any C-hypersurface~$M$ of degree~$d$ in~$\C P^n$ satisfies
$$\chi(M)=\chi^n_d,\quad\sign(M)=\sign^n_d\ .$$
\end{corollary}

{\bf Proof.}

(i) The equality $\chi(M)=\chi^n_d$ follows
immediately by induction from Proposition~\ref{l8}
and the behavior of the Euler characteristic
under ramified coverings.

(ii) Suppose that~$n$ is odd and consider
a nonsingular algebraic hypersurface~$M'$
of degree~$d$ in~$\C P^n$. Let~$m_0$ be as in
Proposition~\ref{l8}. Take any prime number $p > m_0$.
Note that $\widetilde M'=\Pi_p^{-1}(M') \subset \C P^n$
is an algebraic hypersurface of
degree $pd$ tame isotopic to
$\widetilde M = \Pi_p^{-1}(M)$.
The deck transformation
group of the coverings $\Pi_p:\widetilde M\to M$ and
$\Pi_p:\widetilde M'\to M'$
is $G \simeq (\Z/p)^n$.
According to~\cite{Hi2},
\begin{equation}
\sign(M)=\frac{1}{|G|}\sum_{g\in G}\sign(g,\widetilde M),
\quad\sign(M')=\frac{1}{|G|}\sum_{g\in G}\sign(g,\widetilde M')\ .
\label{e11}\end{equation}
For $g=\Id\in G$, we have $\sign(g,\widetilde M)=\sign(\widetilde M)=
\sign(\widetilde M')=\sign(g,\widetilde M')$.
Pick $g\ne\Id$.
The twisted
signatures $\sign(g,\widetilde M)$ and $\sign(g,\widetilde M')$
depend on the embedding of $\Fix(g)$
in~$\widetilde M$ and~$\widetilde M'$, respectively,
and the action of $g$ in the tangent and
normal bundles of $Y=\Fix(g) \subset \widetilde M$
and $Y'=\Fix(g) \subset \widetilde M'$ (see Theorem 6.12 \cite{AS}).
We do not compute the twisted signatures along
Theorem 6.12 \cite{AS},
but use
Hirzebruch's formula for the signature of ramified coverings.
Namely,
let $\langle g\rangle\subset G$ be the cyclic group
of order~$p$ generated by~$g$. According to
\cite{Hi1},
$$p\cdot
\sign(\widetilde M/\langle g\rangle)=\sum_{i=0}^{p-1}\sign(g^i,
\widetilde M)$$
is a universal function of $p$ and the signatures of
self-intersections $Y\circ Y$, $Y\circ Y\circ Y$, $\ldots \ $.
Note that $Y=\Fix(g) \subset \widetilde M$ and
$Y'=\Fix(g) \subset \widetilde M'$
are the intersections of $\widetilde M$ and $\widetilde M'$,
respectively, with the
same
collection of coordinate planes. Hence, due to the tame isotopy of
$\widetilde M$ and $\widetilde M'$, one has
$\sign(Y\circ Y)=\sign(Y'\circ Y')$, $\sign(Y\circ Y\circ Y)=
\sign(Y'\circ Y'\circ Y')$, and so on. Therefore
for any $g\in G$ we have
$$\sum_{i=0}^{p-1}\sign(g^i,
\widetilde M)=\sum_{i=0}^{p-1}\sign(g^i,
\widetilde M')\ ,$$
which immediately implies the
required equality of the right hand sides in
(\ref{e11}), since $G \simeq (\Z/p)^n$ can be decomposed
in the union of cyclic subgroups such that the only intersection
of any two of these subgroups is $\Id \in G$.
\proofend

\begin{corollary}\label{c6}
Let~$M$ be a C-hypersurface of degree~$d$ in~$\C P^n$.
Then
\begin{itemize}
\item $M$ is simply connected if $n>2$,
\item $\pi_1(\C P^n\backslash M)=\Z/d\Z$ if $n\ge 2$.
\end{itemize}
\end{corollary}

{\bf Proof.}

(i) To show that $M$ is simply connected for $n>2$,
note that an algebraic hypersurface of dimension greater than
$1$ is simply connected. Let $m>m_0$
be as in Proposition~\ref{l8}. Consider
a loop $\gam$ in $M$.
We can move it slightly so that it does not meet the
coordinate hyperplanes in $\C P^n$.
Since $\Pi_m^{-1}(M)$ is simply connected,
$m^n\cdot[\gam]$ is contractible in $M$. Similarly,
$(m+1)^n\cdot[\gam]$ is contractible in $M$, and we are done, because
$m^n$ and $(m+1)^n$ are coprime.

(ii) Since an affine nonsingular algebraic hypersurface in
$\C^{n+1}$,
$n\ge 2$, is simply connected, the previous argument shows
that $\hat M \backslash \C P^n \subset\C^{n+1}$ is simply connected,
where~$\hat M$ is a C-hypersurface in~$\C P^{n + 1}$
and $\C P^n$ is a coordinate
hyperplane in $\C P^{n+1}$.

Let a C-hypersurface~$M$ of degree~$d$ in~$\C P^n$ be defined by a
subdivision ${\cal S}:\ T_d^n=\Del_1\cup...\cup\Del_N$ and
collection of numbers ${\cal A}:T_d^n\cap\Z^n\to\C$. Embed $T_d^n$
into $T_d^{n+1}$ as the face $T_d^{n+1}\cap \{i_{n+1}=0\}$, take
the subdivision $\widetilde{\cal S}: \
T_d^{n+1}=\widetilde\Del_1\cup...\cup\widetilde\Del_N$ with
$\widetilde\Del_k$ being the cone over $\Del_k$ with the vertex
$(0,...,0,d)$, $k=1,...,N$, and define $\widetilde{\cal
A}:T_d^{n+1}\cap\Z^{n+1}\to\C$ by $\widetilde
A_{i_1...i_n0}=A_{i_1...i_n}$, $\widetilde A_{0...0d}=-1$,
$\widetilde A_{i_1...i_ni_{n+1}}=0$ if $0<i_{n+1}<d$. These data
define non-degenerate polynomials $\widetilde
F_k(z_1,...,z_{n+1})=F_k(z_1,...,z_n)-z_{n+1}^d$ with Newton
polytopes $\widetilde\Del_k$, $k=1,...,N$, whose charts can be
glued into a C-hypersurface~$\widetilde M$ of degree~$d$ in~$\C
P^{n+1}$. The required isomorphism $\pi_1(\C P^n \backslash
M)=\Z/d\Z$ is a corollary of the following statement.

\begin{lemma}\label{l10}
There is a $\Z/d\Z$-covering $\widetilde M\backslash\{z_{n+1}=0\}\to
\C P^n\backslash M$.
\end{lemma}

{\bf Proof.} We construct a $\Z/d\Z$-covering
$\Phi:
\C Ch(\widetilde{\cal S},\widetilde{\cal A})\backslash\{w_{n+1}=0\}
\to\C T_d^n\backslash\C Ch({\cal S},{\cal A})$ and then show that
it commutes with the $S^1$-action on $\partial
\C Ch(\widetilde{\cal S},\widetilde{\cal A})\backslash\{z_{n+1}=0\}$
and $\partial\C T_d^n\backslash\C Ch({\cal S},{\cal A})$,
thereby defining the required covering.

Let $(w_1,...,w_n,w_{n+1})\in\C Ch(\widetilde{\cal S},
\widetilde{\cal A})$, $w_{n+1}\ne 0$. Then
$(w_1,...,w_{+1})\in\C I(\widetilde\del)$, where $\widetilde\del$
is a cone over a face $\del$ of $\Del_k$, $1\le k\le N$, with
the vertex at $(0,...,0,d)$. Hence $(w_1,...,w_{n+1})=
\tmu_{\widetilde\del}(z_1,...,z_{n+1})$ for some
$(z_1,...,z_{n+1})\in(\C^*)^{n+1}$. So, we define
$$\Phi(w_1,...,w_{n+1})=\tmu_\del(z_1,...,z_n)\ .$$

The map $\Phi$ is well defined. If $\del=\Del_1$, $\ldots$, $\Del_N$,
then $\tmu_{\widetilde\del}$ and $\tmu_\del$ are diffeomorphisms.
If $\del$ lies in an $(n-1)$-plane $\alp_1i_1+...+\alp_ni_n=\bet$,
$i_{n+1}=0$, then $\widetilde\del$ lies in an $n$ plane
$\alp_1i_1+...+\alp_ni_n+\bet i_n/d=\bet$. This means that
$(\tmu_{\widetilde\del})^{-1}(w_1,...,w_{n+1})$ contains the
family $(z_1t^{\alp_1},...,z_nt^{\alp_n},z_{n+1}^{\bet/d})$,
$t\in\R^*_+$, but $\tmu_\del$ takes the family
$(z_1t^{\alp_1},...,z_nt^{\alp_n})$, $t\in\R^*_+$, to one point.

The map $\Phi$ is continuous. If, for some
$k=1,...,N$ and a curve $\widetilde\gam(t)\in(\C^*)^{n+1}$:
$$(\lam_1t^{k_1}+O(t^{k_1+1}),...,\lam_nt^{k_n}+
O(t^{k_n+1}),\lam_{n+1}t^{k_{n+1}}+O(t^{k_{n+1}+1})),
\quad t>0,$$
one has
$\lim_{t\to 0}\tmu_{\widetilde\Del_k}(\widetilde\gam(t))
=\tmu_{\widetilde\del}
(\lam_1,...,\lam_n,\lam_{n+1})$, then
$\lim_{t\to 0}\tmu_{\Del_k}(\gam(t))=\tmu_\del(\lam_1,...,\lam_n)$,
where
$$\gam(t)=(\lam_1t^{k_1}+O(t^{k_1+1}),...,\lam_nt^{k_n}+
O(t^{k_n+1}))\in(\C^*)^n,\quad t>0\ .$$

The map $\Phi$ is surjective. Indeed, if $(w_1,...,w_n)\in
\C T_d^n\backslash\C Ch({\cal S},{\cal A})$, then
$(w_1,...,w_n)=\tmu_\del(z_1,...,z_n)$ for a face
$\del$ of $\Del_k$ ($1\le k\le N$), where $(z_1,...,z_n)\in
(\C^*)^n$, $F^\del_k(z_1,...,z_n)\ne 0$. Then there exists
$z_{n+1}\ne 0$ such that $F^\del_k(z_1,...,z_n)=z_{n+1}^d$;
hence $\tmu_{\widetilde\del}(z_1,...,z_n,z_{n+1})$
belongs to $\C Ch(F_k^{\widetilde\del})
\backslash\{w_{n+1}=0\}$ and $\Phi(\tmu_{\widetilde\del}(z_1,...,
z_n,z_{n+1}))=(w_1,...,w_n)$.

The map $\Phi$ is a $\Z/d\Z$-covering. Indeed, for any point
$(w_1,...,w_n)\in
\C T_d^n\backslash\C Ch({\cal S},{\cal A})$ its preimage
$\Phi^{-1}(w_1,...,w_n)$ consists of $d$ distinct points
$(w'_1,...,w'_n,w'_{n+1}\omega)$ with some fixed $w'_1,...,
w'_{n+1}$ and any $d$-th root of unity $\omega$.

At last, note that if $\Phi(w'_1,...,w'_n,w'_{n+1})=
(w_1,...,w_n)$, then $\Phi(w'_1v,...,w'_nv,w'_{n+1}v)=
(w_1v,...,w_nv)$ for arbitrary $v\in S^1$,
which completes the proof.
\proofend

\begin{proposition}\label{t1}
Let $M \subset \C P^n$ be a C-hypersurface and
$M'\subset \C P^n$
a nonsingular algebraic hypersurface, both of degree
$d$. Then $H_*(M)$ and $H_*(M')$ are isomorphic as graded
groups.
\end{proposition}

{\bf Proof.}
Take~$m_0$ as in Proposition~\ref{l8}. Choose
$m > m_0$ in such a way that
$m$ is coprime with all the orders of elements in
$\Tors H_*(M)$, and
consider $M_m=
\Pi_m^{-1}(M)$.
Pick a homology class~$\alpha \in H_i(M)$.
Note that $m^n\alpha$ belongs to the image
of $(\Pi_m)_*:H_i(M_m) \to H_i(M)$.
For an odd~$i$ different
from $n - 1$, we have $H_i(M_m) = 0$ and, hence,
$H_i(M) = 0$.
In particular, $H_*(M)$ has no torsion.
Any group $H_i(M)$ with even~$i$ ($0 \leq i \leq 2n - 2$) contains
a nontrivial element: the fundamental class of
the intersection of~$M$ with coordinate
hyperplanes taken in appropriate number.
Thus, $H_i(M)$ is isomorphic to~$H_i(M_m) \simeq \Z$ if~$i$
is even and different from $n - 1$.
The fact that the groups $H_{n - 1}(M)$ and $H_{n - 1}(M')$
are isomorphic follows now from the equality
$\chi(M)=\chi(M')$ proven in Corollary \ref{c2}.
\proofend

\begin{proposition}\label{lattice}
Let~$n$ be a positive odd number, $M$ a C-hypersurface of degree~$d$
in~$\C P^n$ and
$M'$ a nonsingular algebraic hypersurface of degree~$d$ in~$\C P^n$.
Then the lattices $(H_{n - 1}(M), B_M)$
and $(H_{n - 1}(M'), B_{M'})$
(where
$B_M: H_{n - 1}(M) \times H_{n - 1}(M) \to \Z$ and
$B_{M'}: H_{n - 1}(M') \times H_{n - 1}(M') \to \Z$ are
the intersection forms on~$M$ and~$M'$, respectively)
are isomorphic.
\end{proposition}

{\bf Proof.}
The lattices
$H_{n-1}(M)$ and $H_{n-1}(M')$ are unimodular and
have the same rank and signature (Corollary \ref{c2}).
It remains to show that these two lattices have
the same parity, since for $d = 1$ we have
$H_{n-1}(M) \simeq H_{n-1}(M') \simeq \Z$, and
for $d \geq 2$ the lattice $H_{n-1}(M')$ is indefinite.
Let us show, first, that
the lattices $H_{n-1}(M)$ and $H_{n-1}(\widetilde M)$ have
the same parity, where $\widetilde M=\Pi_m^{-1}(M)$,
$m_0$ is as in Proposition~\ref{l8} and $m > m_0$ is odd.

(i) If $\alp^2$ is odd, $\alp\in H_{n-1}(M)$,
then $((\Pi_m)^!\alp)^2=m^{2n}\alp^2$ is odd as well.

(ii) If $\bet^2$ is odd,
$\bet\in H_{n-1}(\widetilde M)$, then $\gam^2$ is odd, where
$\gam=\sum_{g\in G}g_*\bet$, and $G=(\Z/m)^n$ is
the deck transformation
group of $\Pi_m$. Furthermore, $\gam=(\Pi_m)^!\alp$ for certain
$\alp\in H_{n-1}(M)$, and $\alp^2$ is odd.

Similarly, the lattices $H_{n-1}(M')$ and $H_{n-1}(\widetilde M')$,
where $\widetilde M' =\Pi_m^{-1}(M')$, have the same parity.
Since $\widetilde M$ and $\widetilde M'$ are isotopic in $\C P^n$,
we obtain that the lattices $H_{n-1}(M)$ and $H_{n-1}(M')$
have the same parity.
\proofend

\begin{corollary}\label{homeo}
A C-surface of degree~$d$ in $\C P^3$ is homeomorphic to
a nonsingular algebraic surface of degree~$d$ in $\C P^3$.
\end{corollary}

{\bf Proof.}
Note that a C-surface of degree~$d$ in~$\C P^3$
is simply connected (see~\ref{c2}),
smoothable (see Proposition
\ref{l7}) and has the same intersection form
as a nonsingular algebraic surface of degree~$d$
in~$\C P^3$ (see~\ref{lattice}).
It remains to apply Freedman's theorem (see, for example,
\cite{FQ}).
\proofend

\subsection{Ramified double coverings}

\begin{proposition}\label{p11}
Let~$M$ be a C-hypersurface of
even degree~$d$ in~$\C P^n$.
Then there exists a closed
simply connected $2n$-dimensional
PL-manifold $Y$ and a map $\Pi:Y\to\C P^n$ which is a double
covering ramified along $M$.
\end{proposition}

{\bf Proof.} The existence of a ramified double covering
$\Pi:Y\to\C P^n$ with $Y$ being a PL-manifold follows generically
from the fact that $[M]\in H_{2n-2}(\C P^n)$ is an even class.
However, we prefer to give an explicit construction of~$Y$.

Namely, we repeat the construction of step (ii) in the proof of
Corollary \ref{c6}, taking
the simplex $T^{n+1}_{d,2} \subset \R^{n+1}$
with vertices $(0,...,0)$,
$(d,0,...,0)$, ..., $(0,...,0,d,0)$, $(0,...,0,2)$ instead of
$T_d^{n+1}$. So, $T_d^n$ embeds into $T_{d,2}^{n+1}$ as the face
$T_{d,2}^{n+1}\cap\{i_{n+1}=0\}$, the subdivision
$T_{d,2}^{n+1}=\widetilde\Del_1\cup...\cup\widetilde\Del_N$ is defined
as the cone over the subdivision $T_d^n=\Del_1\cup...\cup\Del_N$
with the vertex at $(0,...,0,2)$, and $\widetilde{\cal A}:
T_{d,2}^{n+1}\cap\Z^{n+1}\to\C$ is defined by ${\cal A}_{i_1...i_n0}=
A_{i_1...i_n}$, $\widetilde{\cal A}_{i_1...i_n1}=0$,
$\widetilde{\cal A}_{0...02}=-1$. These data define
a PL-manifold $\C Ch(\widetilde{\cal S},\widetilde{\cal A})$
with boundary $\partial\C Ch(\widetilde{\cal S},\widetilde{\cal A})=
\C Ch(\widetilde{\cal S},\widetilde{\cal A})\cap\partial\C
T_{d,2}^{n+1}$, which is the union of the charts of
the polynomials $\widetilde F_k(z_1,...,z_n,z_{n+1})=
F_k(z_1,...,z_n)-z_{n+1}^2$, $k=1,...,N$, and, as in the proof
of Lemma \ref{l10}, there exists a double covering
$\Phi:\C Ch(\widetilde{\cal S},\widetilde{\cal A})\to\C T_d^n$
ramified along $\C Ch({\cal S},{\cal A})$.

Now note that $S^1$ acts on $\partial\C T_{d,2}^{n+1}$ by
$$v\in S^1,\ (w_1,...,w_n,w_{n+1}) \in \partial
\C T_{d,2}^{n+1}\ \mapsto\ (w_1v,...,w_nv,w_{n+1}v^{d/2})
\in \partial\C T_{d,2}^{n+1}\ .$$
The manifold $\partial\C Ch(\widetilde{\cal S},\widetilde{\cal A})$
is invariant with respect to this action, because it
is defined by quasi-homogeneous
polynomials $f(z_1,...,z_{n+1})$ satisfying
$$f(z_1\tau,...,z_n\tau,z_{n+1}\tau^{d/2})=\tau^df(z_1,...,z_{n+1}),$$
which is compatible with the $S^1$-action.
By the same reason
$$\Phi(w_1v,...,w_nv,w_{n+1}v^{d/2})=(w'_1v,...,w'_nv),\quad
v\in S^1\ ,$$
as far as $\Phi(w_1,...,w_n,w_{n+1})=(w'_1,...,w'_m)$,
$(w_1,...,w_{n+1})\in\partial\C Ch(\widetilde{\cal S},
\widetilde{\cal A})$. Hence $\Phi$
reduces to a double covering $\Pi:Y\to\C P^n$, where
$Y=\C Ch(\widetilde{\cal S},\widetilde{\cal A})/S^1$,
ramified along $M$. It remains to show that
$Y$ is
a closed PL-manifold.
Indeed, for any edge of $T_{d,2}^{n+1}$ there
exists a combination of an automorphism of $\Z^{n+1}$,
leaving the $(n+1)$-st axis fixed, and shifts
which puts this edge on a coordinate axis and the adjacent faces
of $T_{d,s}^n$ on the corresponding coordinate planes.
One can easily verify that the orbits of the $S^1$-action on that
edge and adjacent faces contract into points, whose neighborhoods
in $Y$ are
homeomorphic to $\R^{2n}$ as shown in the proof of
Proposition~\ref{l3}.

At last, we show that $\pi_1(Y)=0$. Any loop in $\C P^n$
through a point $\ow\in M$ lifts to a loop on $Y$
through $\ow$, because $\ow$ is covered in $Y$ by itself only.
Since any loop in $\C P^n$ is contractible, so does its lifting
in $Y$.
\proofend

\begin{proposition}\label{p13}
Let~$n$ and~$d$ be positive even numbers,
$M$ be a C-hypersurface
of degree~$d$ in~$\C P^n$,
and~$Y$ be a
double covering of~$\C P^n$ ramified along $M$.
Then
$$\chi(Y)=\chi_{d,2}^n,\quad\sign(Y)=\sign_{d,2}^n,\quad H_*(Y)
\simeq H_*(Y'),$$
and the lattices~$H_n(Y)$ and $H_n(Y')$, equipped with the
intersection forms, are isomorphic, where $Y'$ is the double
covering of~$\C P^n$
ramified along a nonsingular algebraic
hypersurface of degree $d$, and $\chi_{d,2}^n$ and $\sign_{d,2}^n$
are, respectively, the Euler characteristic and
the signature of~$Y'$.
\end{proposition}

{\bf Proof.} The equalities
$\chi(Y)=\chi(Y')$ and $\sign(Y)=\sign(Y')$ follow immediately from
the additivity of the Euler characteristic, the Hirzebruch formula
for the signature of ramified coverings and Corollary \ref{c2}.
The isomorphism $H_*(Y)\simeq H_*(Y')$, by
\cite{Fa}, reduces to
$H_{2i-1}(Y)=0$, $i=1,...,n$, $H_{2i}(Y)=\Z$, $i=0,...,n$, $i\ne n/2$,
$H_n(Y)\simeq\Z^r$, $r=\rk H_n(Y')$.
It can be proven as Proposition \ref{t1}, using a
ramified covering of $Y$ similar to that constructed in Proposition
\ref{l8}. The same ramified covering applied as in the proof
of Proposition \ref{lattice}, gives a lattice isomorphism
$H_n(Y)\simeq H_n(Y')$.
\proofend

\begin{proposition}\label{p12}
Let~$n$ and~$d$ be positive even numbers and~$M$ be
a real C-hypersurface of degree~$d$ in~$\C P^n$.
Then $\R P^n=\R M_+\cup\R M_-$, where $\R M_+$ and $\R M_-$ are
compact manifolds with boundary such that
$$\partial\R M_+=\partial\R M_-=\R M_+\cap\R M_-=\R M\ .$$
The ramified double covering $\Pi:Y\to\C P^n$ admits an action
of the group $\Z/2\Z\oplus\Z/2\Z=\{\Id,\tau,\conj_+,\conj_-\}$,
where $\tau$ is the deck transformation of $\Pi$ and
$\Pi\circ\conj_+=\Pi\circ\conj_-=\conj\circ\Pi$. In addition,
$$\Fix(\tau)=M,\quad\R Y_{\pm}\stackrel{\rm def}{=}
\Fix(\conj_{\pm})=\Pi^{-1}(\R M_{\pm})\ .$$
\end{proposition}

{\bf Proof.}
In the framework of the construction in the proof of
Proposition \ref{p11}, we define $\R Ch({\cal S},{\cal A})_+$
(resp., $\R Ch({\cal S},{\cal A})_-$)
as the union of the
closures of $\tmu_{\Del_k}(\{F_k\ge 0\}\cap(\R^*)^n)$
(resp., $\tmu_{\Del_k}(\{F_k\le 0\}\cap(\R^*)^n)$),
$k=1,...,N$. The action of $S^1$ on
$\partial\C T_d^n$ reduces to the action of $S^0=\{\pm 1\}$ on
$\partial\R T_d^n$ which preserves $\R Ch({\cal S},{\cal A})_+$
and $\R Ch({\cal S},{\cal A})_-$, thus defining
$\R M_{\pm}=\R Ch({\cal S},{\cal A})_{\pm}/S^0$.
The involutions
\begin{eqnarray}
&(w_1,...,w_n,w_{n+1})\mapsto(w_1,...,w_n,-w_{n+1})\ ,\nonumber\\
&(w_1,...,w_n,w_{n+1})\mapsto\conj(w_1,...,w_n,w_{n+1})\ ,\nonumber\\
&(w_1,...,w_n,w_{n+1})\mapsto\conj(w_1,...,w_n,-w_{n+1}) \nonumber
\end{eqnarray}
on the double covering
$\C Ch(\widetilde{\cal S}, \widetilde{\cal A})$ of $\C T^n_d$
ramified along $\C Ch({\cal S}, {\cal A})$ induce
the involutions $\tau$, $\conj_+$ and $\conj_-$, respectively,
on $Y$.
\proofend

\vskip10pt

Now we have to make a digression on
the topology of smooth (not necessarily algebraic)
hypersurfaces in $\R P^n$.
Following V.~Kharlamov~\cite{Kh2} and
O.~Viro~\cite{Vi6}, we define a {\it rank} of
a connected smooth hypersurface
in $\R P^n$ to be the maximal integer~$r$ such that the homomorphism
induced in $r$-dimensional homology with $\Z/2\Z$ coefficients
by the inclusion of the hypersurface into the projective space
is nontrivial.
It is easy to see that the intersection of a hypersurface
of rank~$r$ in $\R P^n$ with a transversal projective subspace~$P$
of dimension $k \geq n - r$
is a hypersurface of rank $r - n + k$ in~$P$.
Similarly, we define a {\it rank} of a connected component
of the complement of a hypersurface in $\R P^n$.
Clearly, the rank of a connected component~$C$ of the complement
of a two-sided hypersurface ({\it i.e.}, a hypersurface
dividing its tubular neighborhood)
is greater or equal to the rank of each component
of the boundary~$\partial C$.

A component~$C$ of the complement of a two-sided hypersurface
is called {\it principal} if the rank of~$C$ is greater then the rank
of each component of~$\partial C$.

\begin{proposition}\label{principal-component}
Let~$S$ be a two-sided hypersurface in $\R P^n$.
Then there exists at most one principal component of
$\R P^n \setminus S$.
\end{proposition}

{\bf Proof.}
Suppose that~$\R P^n \setminus S$ has two principal components~$C_1$
and~$C_2$.
Let~$r_1$ (resp.,~$r_2$) be the maximal rank of a connected component
of~$\partial C_1$ (resp.,~$\partial C_2$).
Assume that $r_1 \leq r_2$.
Consider a projective subspace~$P$ of $\R P^n$ transversal to~$S$
and of dimension~$n - r_1$.
Each connected component of $P \cap \partial C_1$ is of rank~$0$
in~$P$
and divides $P$ into two parts. One of these parts
has rank~$0$ in $P$,
and is called the {\it interior} of the component.
Denote by $C$ the unique component of $P \setminus \partial C_1$
which is not contained  in the interior
of any component of $P \cap \partial C_1$.
Note that $P \cap C_1$ should have a connected component
of positive rank.
Therefore $C \subset P \cap C_1$.
It implies that $P \cap C_2$ is contained in the interior
of some component of $P \cap \partial C_1$,
and thus should have rank~$0$.
\proofend

Clearly, the real parts
of a real algebraic hypersurface of even degree~$d$
and of a real C-hypersurface of even degree~$d$ are two-sided
in $\R P^n$.
Switching if necessary~$\R M_+$ and~$\R M_-$,
we suppose from now on that $\R M_-$ contains
the principal component of $\R P^n \setminus \R M$,
if this principal component does exist.

Denote by~$\defect(\R M)$ the dimension of the intersection
of the kernel of $\frown \omega: H_*(\R M_-, \R M; \Z/2\Z) \to
H_*(\R M_-, \R M; \Z/2\Z)$, defined by $\sigma \mapsto \sigma \frown \omega$,
where
$\omega$ is $(d/2)$-times the generator of $H^1(\R P^n; \Z/2\Z)$,
with the kernel
of the boundary homomorphism
$H_*(\R M_-, \R M; \Z/2\Z) \to H_*(\R M; \Z/2\Z)$.
We call this dimension~$\defect(\R M)$ the {\it defect} of $\R M$.
For any manifold~$X$ denote by~$b_*(X)$ the total Betti number
$\dim_{\Z/2\Z}H_*(X; \Z/2\Z)$ of~$X$.

\begin{proposition}\label{theorem-covering}
Let~$n$ be a positive even number,
and~$M$ a real C-hypersurface of even degree~$d$
in~$\C P^n$.
Then $b_*(\R Y_+) = b_*(\R M)$ and
$b_*(\R Y_-) = b_*(\R M) + 2 \cdot \defect(\R M)$.
\end{proposition}

{\bf Proof.}
The statement can be easily derived
from the Smith exact sequence (see, for example, \cite{Br,Wi})
applied to the deck transformation of $\R Y_{\pm}$
(as it is done in~\cite{Ro1}, \cite{Kh} and~\cite{De-Kh}
in the case of real algebraic
hypersurfaces).
\proofend

\begin{remark}\label{remark-covering}
Note
that if~$n$ and~$d$ are even, then $\R P^n \setminus \R M$ should have
a principal component.
Indeed, if~$\R P^n \setminus \R M$ does not have
a principal component, one has
$b_*(\R Y_+) = b_*(\R M)$ and $b_*(\R Y_-) = b_*(\R M)$,
and thus, $\chi(\R Y_+) \equiv \chi(\R Y_-) \mod 4$.
The last congruence is impossible, because
$\chi(\R Y_{\pm}) = 2\chi(\R M_{\pm})$ and
$\chi(\R M_+) + \chi(\R M_-) = \chi(\R P^n) = 1$.
\end{remark}

\section{Topology of real C-hypersurfaces}

\subsection{Generalized Harnack inequalities}

Let~$M'$ be a nonsingular algebraic hypersurface
of degree~$d$ in $\C P^n$. Denote by~$b_d^n$ the total Betti
number $b_*(M') = \dim_{\Z/2\Z}H_*(M'; \Z/2\Z)$ of $M'$.
One has (see, for example, \cite{D-Kh})
$$b_d^n = \frac{(d - 1)^{n + 1} - (-1)^{n + 1}}{d} + d + (-1)^{n + 1}.$$

\begin{theorem}\label{t2}
For a real C-hypersurface~$M$ of degree $d$ in $\C P^n$, one has
\begin{equation}
b_*(\R M) = b_d^n - 2a(\R M), \label{e19}
\end{equation}
where $a(\R M)$ is a nonnegative integer.
Furthermore, if~$n$ and~$d$ are both even, then
\begin{equation}
b_*(\R Y_-) = b_{d,2}^n - 2a(\R M) + 2(\defect(\R M) - 1)
\leq b_{d,2}^n, \label{e22}
\end{equation}
\begin{equation}
b_*(\R Y_+) = b_{d,2}^n - 2a(\R M) - 2, \label{e22'}
\end{equation}
where $Y$ is the double covering of $\C P^n$ ramified along $M$,
$\R Y_{\pm}$ are
the fixed point sets of two liftings~$\conj_{\pm}$
to~$Y$ of the complex conjugation~$\conj$
in $\C P^n$,
and $b_{d,2}^n$ is the total Betti number of
the double covering of $\C P^n$ ramified along
a nonsingular hypersurface of degree~$d$.
\end{theorem}

{\bf Proof.}
The statement follows from the Smith-Floyd inequality
\begin{equation}
\dim_{\Z/2\Z}H_*(\Fix(\tau); \Z/2\Z) \le
\dim_{\Z/2\Z}H_*(X; \Z/2\Z)\label{e23}
\end{equation}
(where~$X$ is a
compact CW-complex and $\tau:X \to X$ is an involution)
and the congruence
\begin{equation}
 \dim_{\Z/2\Z}H_*(X; \Z/2\Z) \; \equiv \;
\dim_{\Z/2\Z}H_*(\Fix(\tau); \Z/2\Z) \mod 2; \label{e24'}
\end{equation}
see for details \cite{Ro1,Wi}.
We apply the Smith-Floyd inequality
and the above congruence to the involutions
$\conj$ on $M$ and $\conj_{\pm}$ on $Y$.
In order to get~(\ref{e19})
we use Proposition~\ref{t1}.
To prove~(\ref{e22}) and~(\ref{e22'})
we notice
that according to Proposition~\ref{theorem-covering} one has
$b_*(\R Y_+) = b_*(\R M)$ and $b_*(\R Y_-) =
b_*(\R M) + 2\defect(\R M)$.
In addition, since $n$ is even,
$$b_*(Y) = \chi(Y) = 2\chi(\C P^n)-\chi(M)
= 2(n + 1) - (2n - b_*(M)) = 2 + b_*(M)\ ,$$
and it remains to apply Proposition~\ref{p13}.
\proofend

\subsection{Congruences}

\begin{theorem}\label{t4}
Let~$M$ be a real C-hypersurface of degree~$d$ in~$\C P^n$.
If~$n$ is odd then
\begin{equation}
\chi(\R M)\equiv\sign_d^n+\bigg\{\begin{matrix}
&0,\quad\mbox{if}\ a(\R M)=0,\\
&\pm 2,\quad\mbox{if}\ a(\R M)=1\end{matrix}\ \mod 16\ .
\label{e24}\end{equation}
If $n$ and $d$ are even, then
\begin{equation}
\chi(\R M_-)\equiv \frac{\sign_{d,2}^n}2+
\bigg\{\begin{matrix}
&0,\quad\mbox{if}\ a(\R M) = \defect(\R M) - 1,\\
&\pm 1,\quad\mbox{if}\ a(\R M) = \defect(\R M) \end{matrix}\
\mod 8\ ,
\label{e25}
\end{equation}
\begin{equation}
\chi(\R M_+)\equiv \frac{\sign_{d,2}^n}2 \pm 1 \mod 8,
\ \mbox{if}\ a(\R M)=0.
\label{e25'}
\end{equation}
\end{theorem}

\begin{remark}
Note that in the case of even~$n$ and~$d$ the condition
$a(\R M) = \defect(\R M) - 1$ is automatically
fulfilled if $a(\R M) = 0$,
and the condition $a(\R M) = \defect(\R M)$
is automatically fulfilled
if $a(\R M) = 1$.
\end{remark}

{\bf Proof of Theorem~\ref{t4}.}
The statement can be proven as the Rokhlin and
Gudkov-Krahnov-Kharlamov congruences in the algebraic case
\cite{Gu,Kh,Ro1,Wi}.

(i) Let $n = 2k + 1$ and
$a(\R M)=0$. As in \cite{Ro1} (see also \cite{Gu,Wi}),
the latter implies the splitting of $H_{n-1}(M)$ into the
orthogonal sum $H_+\oplus H_-$ of unimodular eigenlattices
corresponding to the eigenvalues $\pm 1$ of
$\conj_*:H_{n-1}(M)\to H_{n-1}(M)$.

\begin{lemma}\label{l11}
The signature of involution
$\sign(\conj_*) = \sign(H_+)-\sign(H_-)$ is equal to
$(-1)^k\chi(\R M)$.
\end{lemma}

{\bf Proof.} By the Atiyah-Singer formula (see \cite{AS,Ro1,Wi})
$\sign(\conj_*)$ is equal to $\R M\circ\R M$, the self-intersection of
$\R M$ in $M$. Let us show that $\R M\circ\R M=(-1)^k\chi(\R M)$.
First, we
smooth~$M$ as in Proposition~\ref{l7}.
Then we take a tangent vector field~$V$ on $\R M_{sm}$ having only
finitely many singular points which are all non-degenerate and lie
outside $\R\del$ for any proper face $\del$ of the polytopes
$\Del_1,...,\Del_N$ in the subdivision of~$T_d^n$. Extend the
vector field~$J(V)$ (where~$J$ is the almost complex structure
on~$M_{sm}$ defined in Proposition~\ref{p9}) to a neighborhood
of~$\R M_{sm}$ in~$M_{sm}$ and slightly move $\R M_{sm}$ along
geodesics in~$M_{sm}$ tangent to the field obtained. The result
has transversal intersection points with $\R M_{sm}$ at the
singular points of~$V$. The intersection indices are equal to the
multiplied by $(-1)^k$ indices of the singular points of~$V$ (see
\cite{Ro1,Wi}). \proofend

Denote by $H_e$ the lattice $H_+$ (resp., $H_-$)
if~$k$ is odd (resp.,~$k$ is even).
From the existence of an almost complex structure on a smoothing
of~$M$ (Proposition~\ref{p9}),
it follows that $H_e$ is even (see, for example, \cite{Wi}).
Let $\sigma_e$ be the signature
of $H_e$.
Since $H_e$ is unimodular and even, the signature
$\sigma_e$ is divisible by $8$.
To finish the proof in the case $a(\R M) = 0$,
it remains to note that according to Lemma~\ref{l11}
we have $\sign(M) - \chi(\R M) = 2\sigma_e$.

The prove
(\ref{e25}) in the case $a(\R M) = \defect(\R M) - 1$, we apply
the same arguments to~$(Y, \conj_-)$
and use the relations
$\chi(\R Y_-) = 2\chi(\R M_-)$
and $b_*(\R Y_-) = b_*(\R M) + 2\defect(\R M)$.

(ii) Let now $n = 2k + 1$ and $a(\R M) = 1$.
In this case (see, for example, \cite{Kh,Wi}),
the discriminants
of~$H_+$ and~$H_-$ are $\pm 2$ (from the Smith theory
we get the inequality $\mid \discr(H_{\pm}) \mid \; \leq 2$ and then
use the congruence $\chi(M) \equiv (-1)^k\sign(M) \mod 4$ to show
that $\mid \discr(H_{\pm}) \mid \; \not= 1$).
Using Lemma~\ref{l11} and the fact that the signature
of an even lattice with discriminant~$\pm 2$ is congruent
to $\pm 1 \mod 8$ we immediately obtain the statement required.

To prove~(\ref{e25}) in the case $a(\R M) = \defect(\R M)$
and~(\ref{e25'})
we again apply the previous arguments
to~$(Y, \conj_-)$ and~$(Y, \conj_+)$.
\proofend

\subsection{Comessatti inequality for real C-surfaces}

\begin{theorem}\label{Comessatti}
Let~$M$ be a real C-hypersurface of degree~$d$ in~$\C P^n$ and
$M'$ be a nonsingular algebraic hypersurface of degree~$d$
in~$\C P^n$.
Then
$$2 - h^{1,1}(M') \leq \chi(\R M) \leq h^{1,1}(M').$$
\end{theorem}

{\bf Proof.}
The arguments are completely similar to the proof of the Comessatti
inequality in the case of
real algebraic surfaces.

Let~$H_+$ and~$H_-$ be again eigenlattices of $H_2(M)$
corresponding to the eigenvalues $\pm 1$ of
$\conj_*:H_2(M) \to H_2(M)$. Denote by~$a^+_{\pm}$ (resp.,
$a^-_{\pm}$) the number of positive (resp., negative) squares
in the diagonal form over~$\Q$
of the restriction of $B: H_2(M) \times H_2(M) \to \Z$ to~$H_{\pm}$.
We have
$$\displaylines{
a^+_+ + a^-_+ + a^+_- + a^-_- = \dim H_2(M), \cr
a^+_+ - a^-_+ + a^+_- - a^-_- = \sign(M), \cr
a^+_+ + a^-_+ - a^+_- - a^-_- = \chi(\R M) - 2, \cr
a^+_+ - a^-_+ - a^+_- + a^-_- = -\chi(\R M).}
$$
The third and fourth equalities follow from the Lefschetz fixed point
theorem and Atiyah-Singer theorem, respectively.
We obtain that
$$4a^-_+ = \dim H_2(M) - \sign(M) + 2\chi(\R M) - 2\geq 0,$$
$$4a^-_- = \dim H_2(M) - \sign(M) - 2\chi(\R M) + 2\geq 0.$$
These inequalities together with the Hodge index relations
give the required statement.
\proofend

\subsection{Topology of real C-curves}

Let~$M$ be an oriented smooth connected closed surface
in~$\C P^2$. Then~$M$ is called a {\it flexible curve of degree~$d$}
(see~\cite{Vi4}) if
\begin{itemize}
\item it realizes $d[\C P^1] \in H_2(\C P^2)$,
\item the genus of~$M$ is equal to $(d - 1)(d - 2)/2$,
\item $M$ is invariant under the complex conjugation,
\item the field of tangent planes to~$M$ on $M \cap \R P^2$
can be equivariantly deformed to the field of lines in $\C P^2$
tangent to $M \cap \R P^2$.
\end{itemize}

According to Propositions~\ref{l6}, \ref{p9} and Corollary~\ref{c2}
(a smoothing of) a real C-curve of degree~$d$ in~$\C P^2$ is
a flexible curve of degree~$d$.
Thus, all the restrictions on the topology of flexible
curves are applicable to real C-curves.
We formulate here in the framework of real C-curves
the principal known restrictions on flexible curves.
An extensive list of restrictions to the topology of flexible curves
can be found in~\cite{Vi4}.

Let us start from definitions. The standard definitions applicable to
real algebraic curves can be naturally extended to real C-curves. A
real C-curve~$A$ of degree~$d$ in~$\C P^2$ is called an {\it
$M$-curve} or {\it maximal} if the real part $\R A$
of $A$ has $(d - 1)(d - 2)/2 + 1$
connected components. A real C-curve~$A$ of degree~$d$ in~$\C P^2$ is
called an {\it $(M - i)$-curve} if~$\R A$ has $(d - 1)(d - 2)/2 + 1 -
i$ connected components. A connected component of the real part
of a real C-curve of degree~$d$ in~$\C P^2$ is called an {\it oval}
if it divides $\R P^2$ into two parts. The part homeomorphic to a
disk is called the {\it interior} of the oval. All the connected
components of the real part of a real C-curve of an even degree
in~$\C P^2$ are ovals. Exactly one connected component of
the real part of a real C-curve of an odd degree
in~$\C P^2$ is not an
oval. This component is called {\it nontrivial}. An oval is {\it
even} (resp., {\it odd}) if it lies inside of an even (resp., odd)
number of other ovals of the curve. The numbers of even and odd ovals of a
curve are denoted by~$p$ and~$n$, respectively. The Euler
characteristic of a connected component of the complement in $\R P^2$
of the real part of a real C-curve is called the
{\it characteristic} of an oval bounding the component from outside.
A component of the complement in $\R P^2$
of the real part of a real C-curve is said to be {\it even}
if each of its inner
bounding ovals contains inside an odd number of ovals.

\begin{theorem}\label{real-C-curves1}

\begin{itemize}
\item {\rm Harnack inequality.}
The number of connected components of the real part
of a real C-curve of degree~$d$ in~$\C P^2$ is at most
$(d - 1)(d - 2)/2 + 1$.
\item {\rm Gudkov-Rokhlin congruence.}
For a maximal real C-curve of degree~$2k$ in~$\C P^2$, one has
$$p - n \equiv k^2 \mod 8.$$
\item {\rm Gudkov-Krahnov-Kharlamov congruence.}
Let~$A$ be a real C-curve of degree~$2k$ in~$\C P^2$.
If~$A$ is an $(M - 1)$-curve, then
$$p - n \equiv k^2 \pm 1 \mod 8.$$
\item {\rm Strengthened Petrovsky inequalities.}
For a real C-curve~$A$ of degree~$2k$ in~$\C P^2$,
one has
$$p - n^- \leq\frac{3k(k - 1)}2 + 1, \quad n - p^- \leq 3k(k - 1),$$
where $p^-$ (resp., $n^-$) is the number of even (resp., odd)
ovals of~$\R A$ with negative characteristic.
\item {\rm Strengthened Arnold inequalities.}
For a real C-curve~$A$ of degree~$2k$ in~$\C P^2$,
one has
$$p^- + p^0 \leq\frac{k^2 - 3k + 3 + (-1)^k}2, \quad
n^- + n^0 \leq\frac{k^2 - 3k + 2}2,$$
where $p^0$ (resp., $n^0$) is the number of even (resp., odd)
ovals of~$\R A$ with characteristic~$0$.
\item {\rm Extremal properties of strengthened Arnold inequalities.}
For a real C-curve of degree~$2k$ in~$\C P^2$,
one has

$p^- = p^+ = 0$, if~$k$ is even and
$p^- + p^0 = (k^2 - 3k + 4)/2$,

$n^- = n^+ = 0$, if~$k$ is odd and
$n^- + n^0 = (k^2 - 3k + 2)/2$.
\end{itemize}
\end{theorem}

A real C-curve~$A$ in~$\C P^2$ is said to be of {\it type I} if its
real part~$\R A$ divides~$A$ into two parts; otherwise, the
curve is of {\it type II}. For a curve~$A$ of type I, the
orientations of two halves of $A \setminus \R A$ induce on~$\R A$ two
opposite orientations which are called {\it complex orientations}.
Note that a real C-curve~$A$ is of type~I if and only if
all the algebraic curves used in the construction of~$A$ are of type~I and
complex orientations on the real parts of these curves can be chosen
in such a way that they induce an orientation of~$\R A$.

A pair of ovals of the real part of a real C-curve in~$\C P^2$
is {\it injective} if one of them is inside of the
other one. A collection of ovals is called a {\it nest} if any two of
them form an injective pair. An injective pair of ovals of a real
C-curve is {\it positive} (resp., {\it negative}) if the complex
orientations of the ovals are induced (resp., are not induced) from
some orientation of the annulus bounded by the ovals. Take an oval of
a real C-curve of type I and of an odd degree in~$\C P^2$, and
consider the M\"obius band which is the complement in~$\R P^2$ of the
interior of the oval. The oval is called {\it positive} (resp., {\it
negative}) if the integer homology class realized in the M\"obius band
by the oval equipped with a complex orientation differs in sign
(resp., coincides) with the class of the doubled nontrivial component
equipped with the complex orientation.

\begin{theorem}\label{real-C-curves2}

\begin{itemize}
\item {\rm Klein congruence.}
Let~$A$ be a real C-curve of type I in~$\C P^2$.
If~$A$ is an ($M - i$)-curve, then
$i \equiv 0 \mod 2$.
\item {\rm Arnold congruence.}
For a real C-curve of type I and of degree~$2k$,
one has
$$p - n \equiv k^2 \mod 4\ .$$
\item {\rm Rokhlin-Mishachev formulae.}
For a real C-curve~$A$ of type I and of degree~$2k$ in~$\C P^2$, one has
$$2(\Pi^+ - \Pi^-) = l - k^2,$$
where $l$ is the number of ovals of~$\R A$,
and~$\Pi^+$ and~$\Pi^-$ are the numbers of positive and negative
injective pairs, respectively.
For a real C-curve~$A$ of type I and of degree~$2k + 1$
in~$\R P^2$, one has
$$2(\Pi^+ - \Pi^-) + \Lambda^+ - \Lambda^- = l - k(k + 1),$$
where~$\Lambda^+$ and~$\Lambda^-$ are the numbers of positive and negative
ovals, respectively.
\item {\rm Kharlamov- Marin congruence.}
Let~$A$ be a real C-curve of degree~$2k$
in~$\C P^2$. If~$A$ is an $(M - 2)$-curve and
$p - n \equiv k^2 + 4 \mod 8$,
then~$A$ is of type I.
\item {\rm Rokhlin inequalities.}
Let~$A$ be a real C-curve of type I and of degree~$2k$
in~$\C P^2$.
If~$k$ is even, then
$4 \nu + p - n \leq 2k^2 - 6k + 8$,
where~$\nu$ is the number of odd nonempty exterior bounding ovals
of even components of $\R P^2 \setminus \R A$.

If~$k$ is odd, then $4 \pi + n - p \leq 2k^2 - 6k + 7$,
where~$\pi$ is the number of odd nonempty exterior bounding ovals
of even components of $\R P^2 \setminus \R A$.
\item {\rm Extremal properties of strengthened Arnold inequalities.}
Let~$A$ be a real C-curve of degree~$2k$ in~$\C P^2$.

If~$k$ is even and
$p^- + p^0 = (k^2 - 3k + 4)/2$, then $A$ is of type I.

If~$k$ is odd and
$n^- + n^0 = (k^2 - 3k + 2)/2$, then $A$ is of type I.
\end{itemize}
\end{theorem}

Harnack inequality in the case of real C-curves
constructed using a primitive ({\it i.e.}, such that all its triangles
are of area~$1/2$) triangulation was, first, proved in~\cite{I1}
and then in a different way in~\cite{Ha}.
Rokhlin-Mishachev formulae in the case of real C-curves
constructed using a primitive triangulation was proved in~\cite{Pa}.

It is interesting that one restriction on real C-curves
is not proved in the case of flexible curves.
This restriction was proved in~\cite{dL1}:

{\it consider a real C-curve~$A$ of degree~$d$ in~$\C P^2$;
if~$A$ is constructed out of threenomials,
then the sum of the depths of any two nests of~$\R A$ is at most~$d/2$.}

In the case of real algebraic curves this statement is known as Hilbert's
theorem and is an immediate corollary of the B\'ezout theorem.

\end{document}